\documentclass[reqno,final,11pt]{amsart}
\usepackage{natbib}  
\usepackage{fancyhdr} 
\usepackage{color} 
\usepackage{hyperref} 
\usepackage{graphicx} 
\usepackage{geometry}
\usepackage[utf8]{inputenc}
\usepackage[T1]{fontenc}
\usepackage{verbatim}

\usepackage{nccmath}
\usepackage{amsthm}
\usepackage{amsmath}
\usepackage{amsfonts}
\usepackage{pifont}
\usepackage{tikz}
\usetikzlibrary{calc,shapes,backgrounds,arrows}
\usepackage{amssymb}
\usepackage{bbm}
\usepackage{mathrsfs}
\usepackage{faktor}
\usepackage[mathscr]{eucal}
\usepackage[main=english,brazilian]{babel}
\usepackage{xcolor} 
\usepackage[final]{pdfpages}

\definecolor{aleacolor}{rgb}{0.16,0.59,0.78}

\hypersetup{
	breaklinks,
	colorlinks=true,
	linkcolor=aleacolor,
	urlcolor=aleacolor,
	citecolor=aleacolor}

\theoremstyle{plain}
\newtheorem{theorem}{Theorem}[section]                                          
\newtheorem{hyp}{Hypothesis}[section]                                          
\newtheorem{proposition}[theorem]{Proposition}                
\newtheorem{lemma}[theorem]{Lemma}
\newtheorem{corollary}[theorem]{Corollary}

\theoremstyle{definition}
\newtheorem{definition}[theorem]{Definition}
\theoremstyle{remark}
\newtheorem{remark}[theorem]{Remark}

\makeatletter \@addtoreset{equation}{section} \makeatother

\begin{document}
	
	\title[Hydrodynamic limit for WASEP with collisions]{Hydrodynamic limit of the multi-component \\slow boundary WASEP with collisions}
	\author{Oslenne Araújo}
	\address{Universidade Federal da Paraíba - UFPB, Department of Mathematics, João Pessoa - PB, 58059-900, Brazil}
	\email{oslenne.araujo@gmail.com} 
	
	\author{Patrícia Gonçalves}
	\address{Instituto Superior Técnico, Department of Mathematics, Av. Rovisco Pais 1, 1049-001, Lisbon}
	\email{pgoncalves@tecnico.ulisboa.pt} 

	\author{Alexandre B. Simas}
	\address{Statistics Program, Computer, Electrical and Mathematical Science and Engineering Division, King Abdullah University of Science and Technology (KAUST)}
	\email{alexandre.simas@kaust.edu.sa} 


	\begin{abstract}
		
		In this article, we study the hydrodynamic limit for a stochastic interacting particle system whose dynamics consists in a superposition of several dynamics:  the exclusion rule, that dictates that no more than one particle per site with a fixed velocity is allowed; a collision dynamics, that dictates that particles at the same site can collide and originate  particles with new velocities such that the linear momentum is conserved; a boundary dynamics that injects and removes particles in the system. This last dynamic destroys the conservation law, and its strength is regulated by a parameter $\theta$. The goal is the derivation of the hydrodynamic limit, and the  boundary conditions  change drastically according to the value of $\theta$.
		
	\end{abstract}
	
	\maketitle
	\section{Introduction}\label{sec:1}
	Stochastic Interacting Particle Systems (SIPS) is an area of probability theory devoted to the mathematical analysis of a collective behavior of continuous-time random walks models subject to constraints. These systems arise in statistical physics, biology, and many other fields of science and were introduced by \cite{Spitzer}.
	A classic problem in this field is to derive the macroscopic laws of thermodynamic quantities of a given physical system, considering microscopic dynamics composed of particles that move according to some prescribed stochastic law. These macroscopic laws are governed by Partial Differential Equations (PDEs) or stochastic PDEs, depending on whether one is looking at the convergence to the mean or at the fluctuations around that mean. Convergence to the mean is a scaling limit, called the \textit{Hydrodynamic Limit} which is described by a solution to a PDE, called the hydrodynamic equation, see \cite{Kipnis/Landim}.
	
\begin{figure}[!htbp]
		\centering
		\begin{tabular}{l}
			\includegraphics[scale=0.58]{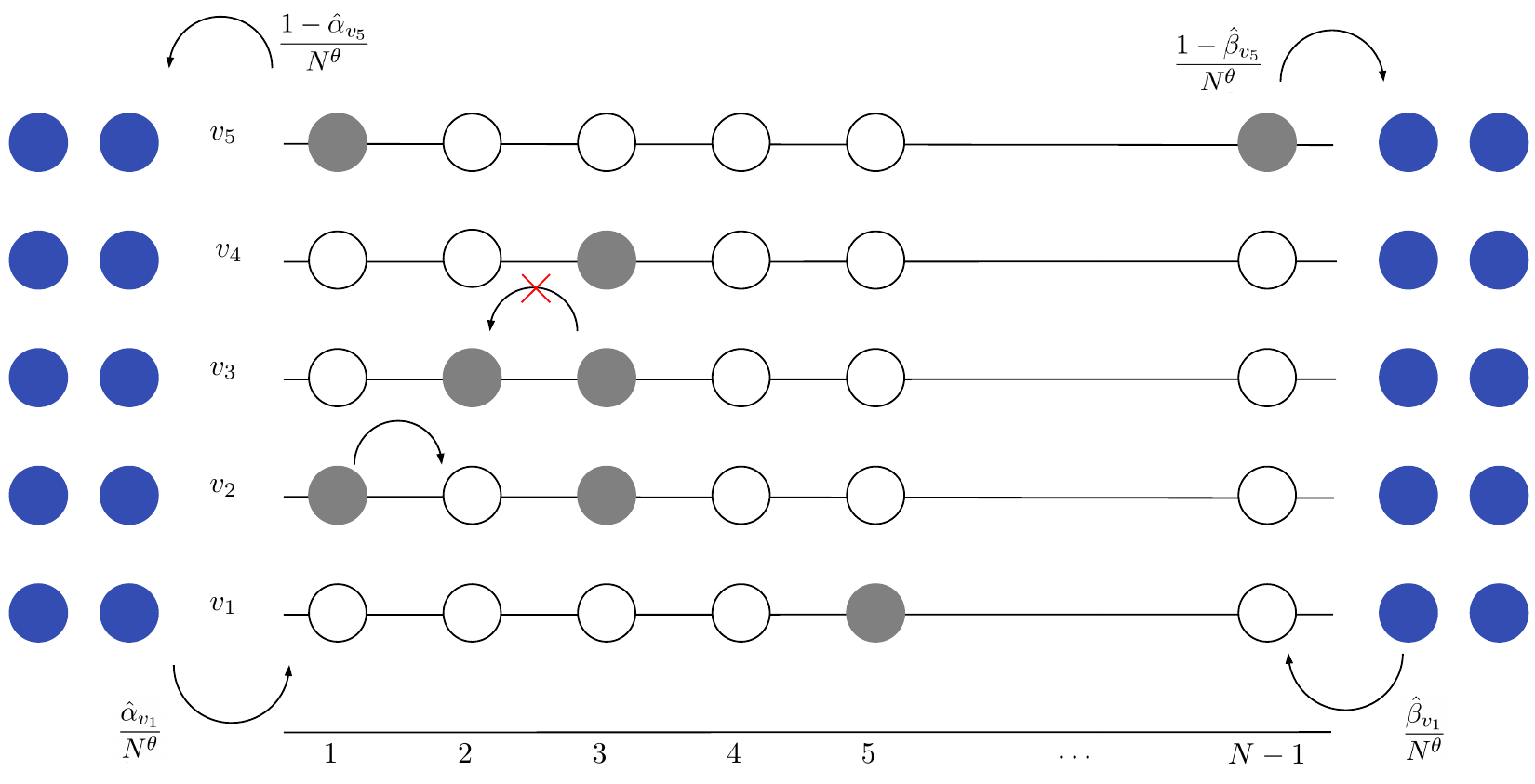} 
		\end{tabular}
		\caption{Illustration of the dynamics}\label{figure1}
	\end{figure}

\begin{figure}[!htbp]
		\begin{center}
			\centering 
			\includegraphics[scale=1.1]{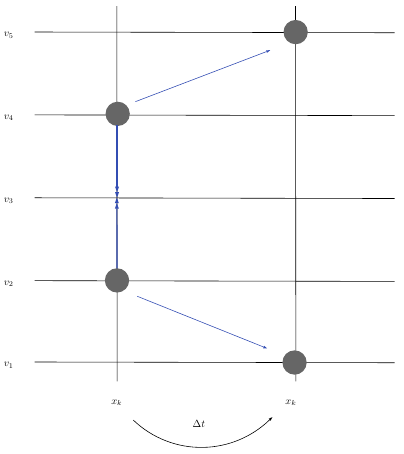} 
			\caption{Illustration of the collision dynamics}\label{figure2}
		\end{center}
	\end{figure}

	To make the reading as enjoyable as possible, we shall informally outline the model that we investigate  in this article when the dimension $d=1$ (see Figures \ref{figure1} and \ref{figure2}), but our results also extend to any dimension $d$.  We consider our particle model evolving on the discrete set of sites $\{1,\dots,N-1\}$, which we call the bulk. Consider now a finite set of possible velocities $\mathcal{V} \subset \mathbb{R}$ and fix a velocity $v \in \mathcal{V}$. At any given time $t$, each site of the bulk is either empty or occupied by one particle with velocity $v$, and  each particle attempts to jump to one of its neighbors with the same velocity, under a weakly asymmetric rate. To prevent the occurrence of more than one particle per site with the same velocity $v$, we introduce an exclusion rule that suppresses  jumps to an already occupied site by a particle with the given velocity $v$. The boundary dynamics is given by the following birth and death process at the sites $1$ or $N-1$:  a particle is inserted into the system at site $1$  with rate $\frac{\hat{\alpha}_v}{N^{\theta}}$  if the site is empty; while if the site $1$ is occupied, a particle is removed from there at rate $\frac{1-\hat{\alpha}_v}{N^{\theta}}$. On the other hand, at site $N-1$ a particle can be  inserted into the system, at rate $\frac{\hat{\beta}_v}{N^{\theta}}$, if the site $N-1$ is empty; while a particle can be removed from the site $N-1$ at the rate $\frac{1-\hat{\beta}_v}{N^{\theta}}$. For the definition of the rates $\hat{\alpha}_v$ and $\hat{\beta}_v$, see \eqref{eq:boundaryrates}. Superposed to these dynamics, there is also a collision process that exchanges velocities of particles at the same site in such a way that the moment is conserved. 
	In Figure \ref{figure1}, we have an illustration of our dynamics: the particles in the bulk are colored in gray, and the particles at the two reservoirs are colored in blue. Note that if a particle at site $x$ with velocity $v$, attempts to jump to an already occupied site $y$ with velocity $v$, the jump is not allowed. In this case, the particle does not move. For example, in Figure \ref{figure1}, the particle at site $3$ with velocity $v_3$ is not allowed to jump to site $2$ since there is already a particle  with velocity $v_3$ at site $2$.
	On the other hand, if the destination site is empty the jump is performed. For example, in Figure \ref{figure1}, the particle at site $1$ with velocity $v_2$ is allowed to jump to site $2$ while keeping the same velocity $v_2$.
	Let us suppose that the clock associated with the leftmost reservoir rings.  Since there exists no particle at site $1$ with velocity $v_1$, a particle can be injected into the system at the site $1$ with velocity $v_1$ at rate $\frac{\hat{\alpha}_{v_1}}{N^{\theta}}$. Also, if the clock associated to the site $1$ with velocity $v_5$ rings, the particle leaves the system at rate $\frac{1-\hat{\alpha}_{v_5}}{N^{\theta}}$ (See Figure \ref{figure1}). Analogously, suppose that the clock associated with the rightmost reservoir rings. Since there is no particle at site $N-1$ with velocity $v_1$, a particle can be injected into the system at the site $N-1$ with the velocity $v_1$ at rate $\frac{\hat{\beta}_{v_1}}{N^{\theta}}$. Also, if the clock associated with the site $N-1$ with velocity $v_5$ rings, the particle leaves the system at rate $\frac{1-\hat{\beta}_{v_5}}{N^{\theta}}$ (See Figure \ref{figure1}).
	Now, let us suppose that the clock associated to $x_k$ rings (see Figure \ref{figure2}). If there are two particles at $x_k$, one with  velocity $v_2$ and another one with velocity $v_4$, then they collide at rate one and produce two particles at the same site $x_k$, but with velocities $v_1$ and $v_5$, satisfying $v_2+v_4=v_1+v_5$, i.e., the preservation of momentum. 
	
	Our goal is to show that the conserved quantities of the system can be described by a hydrodynamic equation. To that end, fix a finite time horizon $[0,T]$, and consider the dynamical behavior of the empirical density and momentum over such interval. The law of large numbers for the empirical density and momentum (the hydrodynamic limit), when the system is taken in the  diffusive scaling limit, is given by a system of parabolic evolution equations of the form
	\begin{equation}\label{eq:PDE}
		\partial_t(\rho,\varrho)+\displaystyle\sum_{v \in \mathcal{V}} \tilde{v}[v \cdot \nabla F(\rho,\varrho)] = \frac{1}{2}\Delta(\rho,\varrho),
	\end{equation}
	where $\tilde{v}=(1, v_1,\dots,v_d)$, $\rho$ stands for the density and $\varrho=(\varrho_1,\dots,\varrho_d)$ stands for the linear momentum. Here, $F$ is a thermodynamical quantity determined by the ergodic properties of the dynamics.
	
	Now we say some words about the proof. We follow the entropy method first introduced in \cite{GPV}, and used in \cite{BPGN}, \cite{Baldasso/Menezes/Neumann/souza} among others. For that, we prove tightness of the sequence of empirical measures in the sense of \cite{mitoma} and then we uniquely characterize the limit point. For the latter, due to the fact that our rates are weakly asymmetric and the system is in contact with slow reservoirs, we need several technical replacement lemmas in order to  write down the notion of weak solutions to each one of the hydrodynamic equations.
	
	The outline of this article is as follows: in Section \ref{sec:3}, we establish the notation adopted in this work and state some useful results and our main theorem, whose proof  is postponed to Section \ref{sec:4} and \ref{sec:5}. In Section \ref{sec:6}, we prove the Replacement Lemmas needed for the proof of the hydrodynamic limit. In Appendix \ref{appendix}, we prove the uniqueness of weak solutions of the hydrodynamic equations, whereas in Appendix \ref{app:diffeo}, we prove an auxiliary result.
	\section{Statement of Results}\label{sec:3}
	We start by introducing some notation to be used in this work. Let $\mathbb{T}^d_N=\{0, \dots,N-1\}^d = (\mathbb{Z}/N\mathbb{Z})^d$ be the $d$-dimensional discrete torus and let $D^d_N$ be the Cartesian product of $S_N$ and $\mathbb{T}^{d-1}_N$, $D^d_N=S_N \times \mathbb{T}^{d-1}_N $, which will henceforth be called \textit{bulk}, with $S_N=\{1, \dots , N-1\}$. Further, let $\mathbb{T}^d = [0,1)^d = (\mathbb{R}/\mathbb{Z})^d$ denote the $d$-dimensional torus and let $D^d=[0,1]\times \mathbb{T}^{d-1}$. 
	
	\begin{hyp}\label{hyp_1}
		Let $\mathcal{V} \subset \mathbb{R}^d$ be a finite set of velocities $v=(v_1, \dots , v_d)$, such that:
		\begin{itemize}\label{eq:conjV}
			\item [(1)] there exists at least one $j\in\{1,\cdots, d\}$ with $v_j\neq 0$;
			\item [(2)] $\mathcal{V}$ is invariant under reflections and permutations of the coordinates, i.e., for any $(v_1, \dots , v_d)\in\mathcal V$ we have that
			$(v_1, \dots , v_{j-1},-v_j,v_{j+1},\dots,v_d)$ and $ (v_{\sigma(1)}, \dots, v_{\sigma(d)})$ belong to $\mathcal{V}$ for all $j\in\{1,\cdots, d\}$, and all permutations $\sigma$ of $\{1, \dots , d\}$.
		\end{itemize}
	\end{hyp}

	The dynamics is chosen in such a way that at most one particle with a certain velocity is allowed at each site of $D^d_N$. We let $\eta(x,v) \in \{0,1\}$ denote the number of particles with velocity $v \in \mathcal{V}$ at site $x \in D^d_N$; $\eta_x=\{\eta(x,v); \, v \in \mathcal{V}\}$ denote the number of particles with velocity $v$ at site $x$; and $\eta=\{\eta_x; \, x \in D^d_N\}$ denotes a configuration. Finally, we let $X_N=(\{0,1\}^{\mathcal{V}})^{D^d_N}$ denote the set of particle configurations.
	
	In the interior of the domain $D^d_N$, the dynamics consists of two parts:
	\begin{enumerate}
		\item particles in the bulk evolve according to nearest-neighbor weakly asymmetric random walks with exclusion among particles with the same velocity, 
		\item binary collisions, that preserve linear momentum, between particles with different velocities. 
	\end{enumerate}
	
	\subsection{The Model}\label{subsec:3.1}
	Our chief goal is to study the stochastic lattice gas model induced by the generator $\mathcal{L}_N$, which is the superposition of the Glauber dynamics with collision among particles of different velocities and an exclusion dynamics (see \cite{Simas} and \cite{Beltran/Landim}):
	\begin{equation}\label{eq:1.1}
		{\mathcal{L}_{N,\theta} = N^2\{\mathcal{L}_{N,\theta}^b+\mathcal{L}_N^c+\mathcal{L}_N^{ex}\},}
	\end{equation}
	where $\theta\geq 0$ is a parameter related to Glauber dynamics. In \eqref{eq:1.1}, $\mathcal{L}^b_{N,\theta}$ denotes the generator of the Glauber dynamics, modelling insertion or removal of particles; $\mathcal{L}_N^c$ denotes the generator that models the collision part of the dynamics; and lastly, $\mathcal{L}_N^{ex}$ models the exclusion part of the dynamics. Note that in \eqref{eq:1.1} time has been speeded up diffusively due to the factor $N^2$.
	
	The jump law and the waiting times are chosen so that the jump rate from site $x$, with velocity $v$, to site $x+y$, with the same velocity $v$, is given by
	\begin{equation}\label{eq:defpn}
		P_N(y,v)=\frac{1}{2}\displaystyle\sum_{j=1}^d(\delta_{y,e_j}+\delta_{y,-e_j})+\frac{1}{N}p(y,v),
	\end{equation}
	where $\delta_{x,y}$ stands for the Kronecker delta, which is equal to one if $x=y$ and $0$ otherwise, and $\{e_1,\dots,e_d\}$ is the canonical basis in $\mathbb{R}^d$. Above  $p(x,v)$ is an irreducible transition probability with finite range and mean velocity $v$, i.e., $\forall v \in \mathcal{V}$, $0 \leq p(x,v) \leq 1, \mbox{~and~} \,\displaystyle\sum_{x \in \mathbb{Z}^d} p(x,v)=1. \mbox{~Moreover,~} \exists \, \mathfrak{K} \in \mathbb{N}:\|x\|>\mathfrak{K}\footnote{$\|\cdot\|$ stands for the Euclidean norm in $\mathbb{R}^d$.}$ implies $p(x,v)=0$, and
	$$
	\displaystyle\sum_{x \in \mathbb{Z}^d}x p(x,v)=v.
	$$
	Note that, since we have finitely many velocities, we can choose a $\mathfrak{K}$ that suits all $v \in \mathcal{V}$. 
	
	The generator of the exclusion part of the dynamics, $\mathcal{L}_N^{ex}$, acts on $f: X_N \to \mathbb{R}$ as
	\begin{equation*}
		(\mathcal{L}_N^{ex}f)(\eta)= \displaystyle\sum_{v \in \mathcal{V}}\sum_{x,z \in D^d_N}\eta(x,v)(1-\eta(z,v))P_N(z-x,v)[f(\eta^{x,z,v})-f(\eta)],
	\end{equation*} 
	where 
	\begin{equation}\label{eq:config}
		\eta^{x,y,v}(z,w)=\left\{
		\begin{array}{ll}
			\eta(y,v), & \mbox{ if } w=v \mbox{ and } z=x,\\
			\eta(x,v), &\mbox{ if } w=v \mbox{ and } z=y,\\
			\eta(z,w), &\mbox{ otherwise}.\\
		\end{array}\right.
	\end{equation}
	
	In view of \eqref{eq:defpn}, we write the decomposition $\mathcal{L}_N^{ex}=\mathcal{L}_N^{ex,1}+\mathcal{L}_N^{ex,2}$ in terms of symmetric and weakly asymmetric part of the generator, where for $f: X_N \to \mathbb{R}$ it holds
	\begin{equation*}\begin{split}
			(\mathcal{L}_N^{ex,1}f)(\eta)=& \displaystyle\frac{1}{2}\sum_{v \in \mathcal{V}}\sum_{\stackrel{x,z \in D^d_N}{|z-x|=1}}\eta(x,v)(1-\eta(z,v))[f(\eta^{x,z,v})-f(\eta)],\\
			(\mathcal{L}_N^{ex,2}f)(\eta)=& \displaystyle\frac{1}{N}\sum_{v \in \mathcal{V}}\sum_{x,z\in D^d_N}\eta(x,v)(1-\eta(z,v))p(z-x,v)[f(\eta^{x,z,v})-f(\eta)].
	\end{split}\end{equation*}
	The generator of the collision part of the dynamics, $\mathcal{L}_N^c$, acts on $f: X_N \to \mathbb{R}$ as 
	\begin{equation*}
		(\mathcal{L}_N^c f)(\eta) = \sum_{y \in D^d_N} \sum_{q \in Q}p_c(y,q,\eta)[f(\eta^{y,q})-f(\eta)],
	\end{equation*}
	where $Q$ is the set of all the possible collisions that preserve linear momentum, i.e.,
	\begin{equation}\label{def:conjQ}
		Q=\{q=(v,w,v',w') \in \mathcal{V}^4 : v+w=v'+w'\}.
	\end{equation}
	The rate $p_c(y,q,\eta)$ is given by
	\begin{equation}\label{eq:defpc}
		p_c(y,q,\eta)=\eta(y,v)\eta(y,w)[1-\eta(y,v')][1-\eta(y,w')],
	\end{equation}
	and for $q=(v_0,v_1,v_2,v_3)$, the configuration $\eta^{y,q}$ after the collision is defined as 
	\begin{equation}\label{eq:configc}
		\eta^{y,q}(z,u)=\left\{ß
		\begin{array}{ll}
			\eta(y,v_{i+2}) &\mbox{ if } z=y \mbox{ and } u=v_i \mbox{ for some } 0 \leq i \leq 3,\\
			\eta(z,u) &\mbox{ otherwise},
		\end{array}\right.
	\end{equation}
	where the index of $v_{i+2}$ should be taken modulo 4.
	Particles of velocities $v$ and $w$ at the same site collide at rate one and produce two particles of velocities $v'$ and $w'$ at the same site but the velocities satisfy the identity $v+w=v'+w'$. 
	
	Finally, the generator of the Glauber dynamics, with parameter $\theta\geq 0$, is given by
	\begin{equation}\label{eq:gen_boundary}
		\begin{split}
			(\mathcal{L}_{N,\theta}^b f)(\eta) & = \displaystyle\frac{1}{2N^\theta}\sum_{\stackrel{x \in D^d_N}{x_1=1}} \sum_{v \in 		\mathcal{V}}\left(\hat{\alpha}_v\left(\tfrac{\tilde{x}}{N}\right)(1-\eta(x,v))+(1-\hat{\alpha}_v\left(\tfrac{\tilde{x}}{N}\right))\eta(x,v)\right)[f(\sigma^{x,v}\eta)-f(\eta)]\\
			& + \displaystyle\frac{1}{2N^\theta}\sum_{\stackrel{x \in D^d_N}{x_1=N-1}} \sum_{v \in \mathcal{V}}\left( \hat{\beta}_v\left(\tfrac{\tilde{x}}{N}\right)(1-\eta(x,v))+(1-\hat{\beta}_v\left(\tfrac{\tilde{x}}{N}\right))\eta(x,v)\right)[f(\sigma^{x,v}\eta)-f(\eta)],
	\end{split}\end{equation}
	where $x=(x_1,\tilde{x})$ with $\tilde{x}= (x_2, \dots, x_d)$,  and
	\begin{equation}\label{eq:configb}
		\sigma^{x,v}\eta(y,w)= \left\{
		\begin{array}{ll}
			1-\eta(x,v), &\mbox{ if } w=v \mbox{ and } y=x,\\
			\eta(y,w), &\mbox{ otherwise, }
		\end{array}
		\right.
	\end{equation}
	for every $v \in \mathcal{V}, \, \hat{\alpha}_v,\, \hat{\beta}_v \in C^2(\mathbb{T}^{d-1})$. We also assume that, for every $v \in \mathcal{V}$, the functions $\hat{\alpha}_v \mbox{ and } \hat{\beta}_v,$ take values in a compact subset of $(0,1)$, and this  implies that $\hat{\alpha}_v(\cdot)$ and $\hat{\beta}_v(\cdot)$ are bounded away from $0$ and $1$. Nevertheless, these conditions can be relaxed in the case $\theta \geq 1$. The functions $\hat{\alpha}_v(\cdot) \mbox{ and }\hat{\beta}_v(\cdot)$ represent the density at the reservoirs. The precise construction of the rates $\hat{\alpha}_v$ and $\hat{\beta}_v$ will be given in the next subsection. The parameter $\theta$ controls the flow of the particles at the boundary of $D^d$. This intuition will be made precise in Theorem \ref{teorema:1}.
	
	In the text, sometimes it will be more convenient to write
	\begin{equation}
		\begin{split}\label{eq:taxas}
			&r^N_x(\eta,\alpha) = \hat{\alpha}_v\left(\tfrac{\tilde{x}}{N}\right)(1-\eta(x,v))+(1-\hat{\alpha}_v\left(\tfrac{\tilde{x}}{N}\right))\eta(x,v), \\
			&r^N_x(\eta,\beta) = \hat{\beta}_v\left(\tfrac{\tilde{x}}{N}\right)(1-\eta(x,v))+(1-\hat{\beta}_v\left(\tfrac{\tilde{x}}{N}\right))\eta(x,v).
		\end{split}
	\end{equation}
	Let $\{\eta_{t,\theta}: \, t \geq 0\}$ be the Markov process with generator $\mathcal{L}_{N,\theta}$ and  denote by $\{S^{N,\theta}_t, \, t \geq 0\}$ the semigroup associated to $\mathcal{L}_{N,\theta}$.
	Let $\mathcal{D}(\mathbb{R}_+, X_N)$ be the set of right continuous functions with left limits taking values in $X_N$ endowed with the uniform topology. For a probability measure $\mu$ on $X_N$, denote by $\mathbb{P}_{\mu,\theta}$ the probability measure on the path space $\mathcal{D}(\mathbb{R}_+, X_N)$ induced by $\{\eta_{t,\theta}: \, t \geq 0 \}$ and the initial measure $\mu$. The expectation with respect to $\mathbb{P}_{\mu,\theta}$ is denoted by $\mathbb{E}_{\mu,\theta}$.
	
	\subsection{Invariant Measures}\label{subsec:3.2}
	In this subsection, we consider the weakly asymmetric exclusion process among particles with the same velocity and collision between particles with different velocities but evolving on the $d$-dimensional discrete torus $\mathbb{T}_N^d$, i.e., with periodic boundary conditions.
	
	For each configuration $\eta \in \{0,1\}^{\mathcal{V}}$ and for $k \in \{0,1,\dots,d\}$, let
	\begin{equation}\label{eq:Ik}
		I_k(\eta)= \sum_{v \in \mathcal{V}} v_k\eta(v), 
	\end{equation}
	where we use the convention that $v_0=1$ so that $I_0(\eta)$ denotes the mass of $\eta$ and for $k \in \{1, \dots , d\}$ $I_k(\eta)$  denotes the $k$-th momentum of $\eta$. Set $\textbf{I}(\eta):=(I_0(\eta), \dots, I_d(\eta))$. Fix $L \geq 1$ and a configuration $\eta$. Let $\textbf{I}^L(x, \eta)=(I^L_0(x), \dots , I^L_d(x))$ be the average of the conserved quantities in a cube of length $L$ centered at $x$:
	\begin{equation}\label{eq:def_boldI}
		\textbf{I}^L(x, \eta)=\frac{1}{|\Lambda_L|}\sum_{z \in x+\Lambda_L}\textbf{I}(\eta_z),
	\end{equation}
	where $\Lambda_L=\{-L, \dots,L\}^d$ and $|\Lambda_L|=(2L+1)^d$ is the discrete volume of the box $\Lambda_L$.
	
	Assume that the set of velocities is chosen in such a way that the unique conserved quantities by the random walk dynamics described above are the mass and the momentum: $\sum_{x \in D^d_N}\textbf{I}(\eta_x)$.
	The following examples of sets of velocities satisfying these conditions were introduced in \cite{Esposito}:
	\begin{enumerate}
		\item[] {\bf Model I:} Denote by $\mathcal{E} = \{e=\pm e_j \mbox{ for some } j \in \{1,\dots,d\}\}$ and let $\mathcal{V}=\mathcal{E}$. Under these choices, the only possible collisions are $q=(v,w,v',w')$ such that $v+w=0$ and $v'+w'=0$. So, we have $\mathcal{V} = \{\pm e_1,\ldots, \pm e_d\}$.\\
		\item[] {\bf Model II:} Let $d=3$ and denote by $\sigma$ any element of $\mathcal{J}$, the permutation group of $\{1,2,3\}$. Let $\mathcal{V}=\{v : \sigma v = \varpi\mbox{~for some~} \sigma \in \mathcal{J}\} $, where $\varpi$ is the positive solution of $\varpi^4-6\varpi^2-1=0$.
	\end{enumerate}
	
	\quad
	
	The proof that the only collision invariants are the total mass and momentum is part of the ergodic theorem also proved in \cite{Esposito}. 
	
	For each chemical potential $\lambda= (\lambda_0, \dots, \lambda_d) \in \mathbb{R}^{d+1}$, let $m_{\lambda}$ denote the probability measure  given on $\xi\in\{0,1\}^{\mathcal{V}}$ by
	\begin{equation}\label{eq:1.2}
		m_{\lambda}(\xi)= \frac{1}{Z(\lambda)}\exp\{\lambda \cdot\textbf{I}(\xi)\}, 
	\end{equation}
	where $Z(\lambda)$ is a normalizing constant. Note that $m_{\lambda}$ is a product measure on $\{0,1\}^{\mathcal{V}}$, i.e., the variables $\{\eta(v): \, v \in \mathcal{V}\}$ are independent under $m_{\lambda}$.
	
	Let $\mu_{\lambda}^N$ denote the product measure on $X_N$ with marginals given by 
	$$
	\mu_{\lambda}^N\{\eta: \, \eta(x, \cdot)=\xi \}= m_{\lambda}(\xi),
	$$ 
	for each $\xi \in \{0,1\}^\mathcal{V}$ and $x \in D^d_N$. Note that $\{\eta(x,v): x \in D^d_N, \, v \in \mathcal{V}\}$ are independent random variables under $\mu_{\lambda}^N$, and that the measure $\mu_{\lambda}^N$ is invariant for the  process when taken with periodic boundary conditions, i.e., with generator given by $\mathcal{L}_N = N^2\{\mathcal{L}_N^c+\mathcal{L}_N^{ex}\}$. The expectation under $\mu_{\lambda}^N$ of the mass and $k$-th momentum are, respectively, given by
	$$
	\begin{array}{l}
		\rho(\lambda):= E_{\mu_{\lambda}^N}[I_0(\eta_x)]=\displaystyle\sum_{v \in \mathcal{V}}\theta_v(\lambda)\quad\textrm{and}\quad
		\varrho_j(\lambda):= E_{\mu^N_{\lambda}}[I_j(\eta_x)]=\displaystyle\sum_{v \in \mathcal{V}} v_j\theta_v(\lambda), \, j\in\{1,\dots,d\}.
	\end{array}
	$$
	In the last formula, $\theta_v(\lambda)$ is the expected value of the density of particles with velocity $v$ under $m_{\lambda}$:
	\small\begin{equation}\label{theta}
		\theta_v(\lambda):= E_{m_{\lambda}}[\xi(v)]= \displaystyle\frac{\exp\left\{\lambda_0+\displaystyle\sum_{j=1}^{d}\lambda_jv_j\right\}}{1+\exp\left\{\lambda_0+\displaystyle\sum_{j=1}^d \lambda_jv_j\right\}}.
	\end{equation}\normalsize
	Let $(\rho,\varrho)(\lambda):= (\rho(\lambda),\varrho_1(\lambda), \dots, \varrho_d(\lambda))$ be the map that associates the chemical potential to the vector of density and momentum. We prove, in Appendix \ref{app:diffeo}, that under Hypothesis \ref{eq:conjV}, the map $(\rho,\varrho)$ is a diffeomorphism onto $\mathfrak{U} \subset \mathbb{R}^{d+1}$, the interior of the convex envelope of $\{I(\xi), \, \xi \in \{0,1\}^{\mathcal{V}}\}$. Denote by $\Lambda=(\Lambda_0, \dots , \Lambda_d): \mathfrak{U} \to \mathbb{R}^{d+1}$ the inverse of $(\rho,\varrho)$. This correspondence allows one to parameterize the invariant states by the density and momenta: for each $(\rho,\varrho) \in \mathfrak{U}$, we have a product measure $\nu_{\rho,\varrho}^N = \mu_{\Lambda(\rho,\varrho)}^N$ on $X_N$.
	
	It is possible to show, after a long and tedious computation, that for \emph{Model I} and for every $v\in\mathcal{V}$, the function $\Phi_v:\mathfrak{U}\to (0,1)$ given by $\Phi_v(\rho,\varrho) = \chi(\theta_v(\Lambda(\rho,\varrho)))$ is Lipschitz, where $\chi(a)=a(1-a)$ is the static compressibility. The details can be found, e.g., in \cite{Farfan/Simas/Valentim}.
	
	\begin{hyp}\label{hyp_2}
		In what follows we will always assume that our set of velocities $\mathcal V$ is such that: for any $v\in\mathcal{V}$ the function $\Phi_v(\cdot)$ is Lipschitz. 
	\end{hyp}
	\begin{remark}
		Note that the previous discussion is restricted to considering the model evolving on the torus. For our model, due to the Glauber dynamics, these product measures are no longer invariant. Nevertheless, we introduced them and recalled their definitions since they will play a key role in the proof of the replacement lemmas, since they will be used as reference measures that are close to the (unknown) invariant measure of the system (see  Section \ref{sec:6}).
	\end{remark}
	
	We are now in a position to define $\hat{\alpha}_v(\cdot)$ and $\hat{\beta}_v(\cdot)$, $v\in\mathcal{V}$, properly. Indeed, let, for each $v\in\mathcal{V}$, $\alpha_v(\cdot),\beta_v(\cdot) \in C^2(\mathbb{T}^{d-1})$, and
	$$
	\alpha(\tilde{u}) := \displaystyle\sum_{v \in \mathcal{V}}\tilde{v}\alpha_v(\tilde{u})\quad \hbox{ and }\quad \beta(\tilde{u}) := \displaystyle\sum_{v \in \mathcal{V}}\tilde{v}\beta_v(\tilde{u}),
	$$
	where $\tilde{u}\in\mathbb{T}^{d-1}$, \, $\tilde{v}=(v_0,v_1,\dots,v_d)$ and  $v_0$  was defined right after \eqref{eq:Ik}.
	Then, we define
	\begin{equation}\label{eq:boundaryrates}
		\hat{\alpha}_v(\tilde{u}) := \theta_v(\Lambda(\alpha(\tilde{u})))\hbox{ and }\hat{\beta}_v(\tilde{u}) := \theta_v(\Lambda(\beta(\tilde{u}))), \quad \tilde{u}\in\mathbb{T}^{d-1}, v\in\mathcal{V}.
	\end{equation}
	Observe that, by definition of $\Lambda(\cdot)$,
	\begin{equation}\label{eq:sumalphabetahat}
		\begin{split}
			&\sum_{v\in\mathcal{V}} \hat{\alpha}_v(\tilde{u}) = \sum_{v\in\mathcal{V}} \alpha_v(\tilde{u}), \quad \sum_{v\in\mathcal{V}} v_j \hat{\alpha}_v(\tilde{u}) = \sum_{v\in\mathcal{V}} v_j \alpha_v(\tilde{u}),\\
			&\sum_{v\in\mathcal{V}} \hat{\beta}_v(\tilde{u}) = \sum_{v\in\mathcal{V}} \beta_v(\tilde{u}) \hbox{ and } \sum_{v\in\mathcal{V}} v_j \hat{\beta}_v(\tilde{u}) = \sum_{v\in\mathcal{V}} v_j \beta_v(\tilde{u}), \tilde{u}\in\mathbb{T}^{d-1},\, j=1,\dots,d.
		\end{split}
	\end{equation}
	
	\subsection{Hydrodynamic Equations}\label{subsec:3.3}
	Recall that we fixed a finite time horizon $[0,T]$. Let $C^{m,n}([0,T]\times D^d)$ be the set of functions defined on $[0,T]\times D^d$ that are $m$ times continuously differentiable on the first variable (the time variable), and $n$ times continuously differentiable on the second variable (the space variable). For a function $G:=G(t,u) \in C^{m,n}([0,T]\times D^d)$, we let $\partial_tG$ denote its derivative with respect to the time variable $t$ and $\partial_{u_i}G$ denote its derivative with respect to the space variable $u_i$, $i\in\{1,\dots,d\}$. For simplicity of notation, we set $\Delta G:=\sum_{i=1}^d \frac{\partial^2 G}{\partial{{u_i}^2}}$ and $\nabla G = (\partial_{u_1},\dots,\partial_{u_d})$ represents the Laplacian operator and generalized gradient of the function $G$, respectively. Finally, $C^{m,n}_0([0,T]\times D^d)$ denotes the set of functions $G \in C^{m,n}([0,T]\times D^d)$ such that for any time $t$ the function $G_t$ vanishes at the boundary of $D^d$, that is, $G_t(0,\tilde{u})=G_t(1,\tilde{u})=0$, where $u \in D^d$ is decomposed as $(u_1,\tilde{u})$, with $\tilde{u} \in \mathbb{T}^{d-1}$.
	
	Let $(B,\|\cdot\|_B)$ be a separable Banach space. We let $L^2([0,T],B)$ be the Banach space of measurable functions $U:[0,T] \to B$ for which
	$
	\|U\|^2_{L^2([0,T],B)}=\displaystyle\int_0^T\|U_t\|^2_B\, dt <\infty.
	$
	Moreover, let $\mathscr{H}^1(D^d)$ be the Sobolev space of measurable functions in $L^2(D^d)$ that have generalized derivatives in $L^2(D^d)$.
	
	We are now in a position to define the system of partial differential equations with different boundary conditions, and the respective notions of weak solutions that appear  in the hydrodynamic limit of our model. We begin by introducing the hydrodynamic equation:
	\begin{definition}\label{defi:hydrod_eq}
		Let $\mathfrak C\in\{0,1,+\infty\}.$ Fix measurable functions $\rho_0: D^d \to \mathbb{R}_+$, $\varrho_0: D^d \to \mathbb{R}^d$ such that {$(\rho_0,\varrho_0)(u) \in \mathfrak{U}$}. The following system of parabolic partial differential equations is  the hydrodynamic equation associated to the stochastic lattice gas model with generator \eqref{eq:1.1} and boundary rates \eqref{eq:boundaryrates}:
		\small\begin{equation}\label{eq:1.6}
			\left\{
			\begin{split}
				&	\partial_t(\rho,\varrho)+\displaystyle\sum_{v \in \mathcal{V}} \tilde{v}[v \cdot \nabla\chi(\theta_v(\Lambda(\rho,\varrho)))] = \frac{1}{2}\Delta(\rho,\varrho),\\
				&	\Big(\displaystyle\frac{\partial(\rho,\varrho)}{\partial u_1}(t,0, \tilde{u}) -2\displaystyle\sum_{v \in \mathcal{V}}\tilde{v}[v_1\, \chi(\theta_v(\Lambda(\rho,\varrho)))]\Big)=\mathfrak C\Big((\rho,\varrho)(t,0,\tilde{u})-\alpha(\tilde{u})\Big), \\
				&	\Big(\displaystyle\frac{\partial(\rho,\varrho)}{\partial u_1}(t,1, \tilde{u})-2\displaystyle\sum_{v \in \mathcal{V}}\tilde{v}[v_1\, \chi(\theta_v(\Lambda(\rho,\varrho)))]\Big)=\mathfrak C\Big(\beta(\tilde{u})-(\rho,\varrho)(t,1,\tilde{u})\Big), \\
				&	(\rho,\varrho)(0, u)=(\rho_0, \varrho_0)(u),
			\end{split}
			\right.
		\end{equation}\normalsize
		where $t \in (0,T],\tilde{u} \in \mathbb{T}^{d-1}$, $u\in D^d$, $\chi(a)=a(1-a)$ is the static compressibility of the system and for each velocity $v=(v_1,\dots,v_d)$, we denote $\tilde{v}=(v_0,v_1,\dots,v_d), \hbox{ with } v_0=1$. 
	\end{definition}
	
	We will now give  the precise notion of weak solution to the system \eqref{eq:1.6}:
	
	\begin{definition}\label{defi:2}
		Let $\mathfrak C\in\{0,1,+\infty\}.$ Fix measurable functions $\rho_0: D^d \to \mathbb{R}_+$, $\varrho_0: D^d \to \mathbb{R}^d$ such that {$(\rho_0,\varrho_0)(u) \in \mathfrak{U}$}. We say that $(\rho, \varrho ): [0,T] \times D^d \to \mathbb{R}_+ \times \mathfrak{U}$ is a weak solution of the system of parabolic partial differential equations \eqref{eq:1.6} if the following two conditions hold:
		\begin{enumerate}
			\item $(\rho, \varrho ) \in L^2([0,T],\mathscr{H}^1(D^d));$
			\item $(\rho,\varrho)$ satisfies the weak formulation:
		\end{enumerate}
		\begin{equation}\begin{split}\label{eq:1.7}
				\mathfrak{F}_{\mathfrak{C}}((\rho,\varrho),t,G)&:=\displaystyle\int_{D^d}(\rho,\varrho)(T,u) G(T,u) \,du-\displaystyle\int_{D^d}(\rho,\varrho)(0,u) G(0,u) \,du 
				\\ 
				&-\displaystyle\int_0^T \int_{D^d}(\rho,\varrho)(t,u)\left( \partial_tG(t,u)+\frac{1}{2}\sum_{i=1}^d \frac{\partial^2G}{\partial u_i^2}(t,u)\right) \,du\,dt \\ 
				&-\displaystyle\int_0^T \int_{D^d}\displaystyle\sum_{v \in \mathcal{V}}\tilde{v}\, \chi(\theta_v(\Lambda(\rho,\varrho)))\sum_{i=1}^d v_i \frac{\partial G}{\partial u_i}(t,u) \, du\, dt\\ 
				&-\displaystyle\frac{\tilde{C}(\mathfrak C)}{2}\displaystyle\int_0^T\!\!\!\int_{\{1\}\times \mathbb{T}^{d-1}}\left(\beta(\tilde{u})-(\rho,\varrho)(t,1,\tilde{u})\right)G(t,1,\tilde{u})\, dS\, dt\\ 
				&+\displaystyle\frac{\tilde{C}(\mathfrak C)}{2}\displaystyle\int_0^T\!\!\!\int_{\{0\}\times \mathbb{T}^{d-1}}\left((\rho,\varrho)(t,0,\tilde{u})-\alpha(\tilde{u})\right)G(t,0,\tilde{u})\, dS\, dt\\ 
				&+\displaystyle\frac{1}{2}\displaystyle\int_0^T\int_{\{1\}\times \mathbb{T}^{d-1}}\Big\{(\rho,\varrho)(t,1,\tilde{u}) \mathbbm{1}_{\{\mathfrak C=0,1\}}+\beta(\tilde u)\mathbbm{1}_{\{\mathfrak C=+\infty\}}\Big\}\frac{\partial G}{\partial u_1}(t,1,\tilde{u})\, dS\, dt \\
				&-\displaystyle\frac{1}{2}\displaystyle\int_0^T\int_{\{0\}\times \mathbb{T}^{d-1}}\Big\{(\rho,\varrho)(t,0,\tilde{u}) \mathbbm{1}_{\{\mathfrak C=0,1\}}+\alpha(\tilde u)\mathbbm{1}_{\{\mathfrak C=+\infty\}}\Big\}\frac{\partial G}{\partial u_1}(t,0,\tilde{u})\, dS\, dt =0
		\end{split}\end{equation}\normalsize
		for all $t \in [0,T]$ and any function $G: [0,T]\times D^d \to \mathbb{R}^{d+1}$ in $C(\mathfrak C)$, where
		$$C(\mathfrak C):=\mathbbm{1}_{\{0,1\}}(\mathfrak C) C^{1,2}([0,T] \times D^d)+\mathbbm{1}_{\{+\infty\}}(\mathfrak C) C_0^{1,2}([0,T] \times D^d)
		\mbox{~and~} \tilde{C}(\mathfrak{C})=\left\{\begin{array}{l}
			0, \mbox{~if~} \mathfrak{C}=+\infty,\\
			1, \mbox{~if~} \mathfrak{C}=1,\\
			0, \mbox{~if~} \mathfrak{C}=0.
		\end{array}
		\right.$$
		Above, for  $u=(u_1,\tilde{u}) \in \{0,1\} \times \mathbb{T}^{d-1}$, and $S$ denotes the Lebesgue measure on $\mathbb{T}^{d-1}\cong [0,1)^{d-1}$.
	\end{definition}
	\begin{remark}\label{obs:1}
		In Definition \ref{defi:2}, if $\mathfrak C=0$, we say that $(\rho, \varrho )$ has Neumann boundary conditions, if $\mathfrak C=1$, we say that $(\rho, \varrho )$ has Robin boundary conditions, and if $\mathfrak C=+\infty$, we say that $(\rho, \varrho )$ has Dirichlet boundary conditions that we interpret as 
		$ (\rho,\varrho)(t,u)=\left\{\begin{array}{l}
			\alpha(\tilde{u}), \mbox{~if~} u_1=0,\\
			\beta(\tilde{u}), \mbox{~if~} u_1=1.
		\end{array}\right.$ We note however that these equalities at the boundary, in the Dirichlet case,  are in the sense of the trace of the operator, since condition (1) in Definition \ref{defi:2} above is in force.
	\end{remark}
	
	Our proof boils down to applying carefully the  entropy method  developed by \cite{GPV}. To that end we need to prove the uniqueness of the weak solutions of our hydrodynamic equations. We could not cover all the cases but  we refer the reader to Appendix \ref{appendix} for a proof.
	
	\subsection{Hydrodynamic Limit}\label{subsec:3.4}
	Let $\mathcal{M}_+$ be the space of finite positive measures on $D^d$ endowed with the weak topology and let $\mathcal{M}$ be the space of bounded variation signed measures on $D^d$ endowed with the weak topology. Let $\mathcal{M}_+\times \mathcal{M}^d$ be the Cartesian product of these spaces endowed with the product topology, which is metrizable, since the space of continuous functions defined on $D^d$, denoted by $C(D^d)$ is separable, and by the Riesz-Markov theorem, we have $C(D^d)^{*}=\mathcal{M}$.
	
	Let $\theta\geq 0$ and recall that $\{\eta_{t,\theta}\}_{t\geq 0}$ is the Markov process with generator $\mathcal{L}_{N,\theta}$. To simplify the notation, for $k\in\{0, \dots, d\}$ we consider $\varrho^0=\rho$ and $\varrho=(\varrho^1,\dots,\varrho^d)$ and let $\pi^{k,N}_{t,\theta}$ denote the empirical measure associated to the $k$-th quantity of interest:
	\begin{equation}\label{eq:1.11}
		\pi^{k,N}_{t,\theta}(du)  =\displaystyle\frac{1}{N^d}\sum_{x \in D^d_N}I_k(\eta_{t,\theta}(x))\delta_{\frac{x}{N}}(du),
	\end{equation}
	where $\delta_u (du)$ stands for the Dirac measure concentrated at $u \in D^d$. We let $\langle \pi^{k,N}_{t,\theta}, G \rangle $ denote the integral of a test function $G$ with respect to the empirical measure $\pi^{k,N}_{t,\theta}$, and for a measure $\nu$ in $X_N$ and $f,g \in L^2(D^d,\nu)$, the symbol $\langle f, g \rangle_{\nu}$ denotes the inner product in $L^2(D^d,\nu)$ i.e.  $\langle f, g \rangle_{\nu} = \int_{D^d} fg \, d\nu.$
	
	Let $\mathcal{D}([0,T], \, \mathcal{M}_+ \times \mathcal{M}^d )$ be the set of right continuous functions with left limits taking values in $\mathcal{M}_+ \times \mathcal{M}^d$ and endowed with the uniform topology. We consider the sequence $(\mathbb{Q}_{N,\theta})_{N\geq1}$ of probability measures on $\mathcal{D}([0,T], \, \mathcal{M}_+ \times \mathcal{M}^d )$ that corresponds to the Markov process $\pi^N_{t,\theta}=(\pi^{0,N}_{t,\theta}, \dots , \pi^{d,N}_{t,\theta})$ starting from a probability measure $\mu^N$. At this point we need to fix initial measurable profiles $\rho_0:D^d \to \mathbb{R}_+$ and $\varrho_0:D^d \to \mathbb{R}^d$, where $\varrho_0=(\varrho_{0,1}, \dots,\varrho_{0,d})$ and an initial distribution $(\mu^N)_N$ on $X_N$. Before introducing the main result, let us define some notions that will play an important role in what follows.
	\begin{definition}[Finite energy]\label{defi:finiteenergy}
		We say that $(\rho, \varrho)\in L^2([0,T],(L^2(D^d))^{d+1})$ has finite energy if its components belong to $L^2([0,T], \mathscr{H}^1(D^d))$, i.e., if $\nabla\rho$ and $\nabla \varrho$ are measurable functions such that for $k \in \{0,1,\dots, d\}$ 
		\begin{equation}\nonumber
			\displaystyle\int_0^T ds \left(\int_{D^d} \| \nabla \varrho^k(s,u)\|^2 du\right)< \infty.
		\end{equation}
	\end{definition}
	
	\begin{definition}\label{defi:4}
		We say that a sequence of probability measures $(\mu^N)_N$ on $X_N$ is associated to the density profile $\rho_0$ and to the momentum profile $\varrho_0$, where we have $Im((\rho_0,\varrho_0)) \subset \mathfrak{U}$, if, for every $G\in C(D^d)$ and for every $\delta >0$,
		\begin{equation}\nonumber
			\displaystyle\lim_{N \to \infty} \mu^N\left[\eta :\left|\frac{1}{N^d}\sum_{x \in D^d_N}G\left(\tfrac{x}{N}\right)I_k(\eta(x))-\int_{D^d}G(u)\varrho_0^k(u)du\right|>\delta \right]=0,
		\end{equation}
		for every $k \in \{0,1,\dots,d\}$.
	\end{definition}
	
	\begin{theorem}\label{teorema:1}
		Fix a set of velocity $\mathcal{V}$ satisfying Hypotheses \ref{hyp_1} and \ref{hyp_2}. Let $\rho_0$ and $\varrho_0$ be measurable functions, also let $(\mu^N)_N$ be a sequence of probability measures on $X_N$ associated to the profile $(\rho_0, \varrho_0)$. Then, for every $t \in [0,T]$, for every function $G \in C(D^d)$, and for every $\delta > 0$,
		$$
		\lim_{N \to \infty} \mathbb{P}_{\mu^N}\left[\eta_. \in \mathcal{D}([0,T], X_N):\left|\displaystyle\frac{1}{N^d}\sum_{x \in D^d_N} G\left(\tfrac{x}{N}\right)I_k(\eta_{t,\theta}(x))- \int_{D^d}G(u)\varrho^k(t,u)\, du\right|> \delta \right]=0,
		$$
		for $k \in \{0,1,\dots,d\}$ where $(\rho,\varrho)$ has finite energy (see Definition \ref{defi:finiteenergy}) and it is the unique weak solution of \eqref{eq:1.6} as given in Definition \ref{defi:2} with  $\mathfrak C=+\infty$, if $0\leq \theta < 1$; $\mathfrak C=1$, if $\theta =1$ and $d=1$; or $\mathfrak C=0$, if $\theta >1$.
	\end{theorem}
	
	\begin{remark}
		In the case $\theta<0$, one can use the same arguments as those in the proof of Theorem \ref{teorema:1}, but taking test functions in $C^{1,2}_c([0,T]\times D^d)$, that is, functions of class $C^1$ in time, and of class $C^2$ in space and compactly supported. Nevertheless, since the space of test functions is too small we need to ask from the notion of weak solution an additional point (3) which, for uniqueness,  requires the solution to satisfy the Dirichlet conditions, almost surely in time. This can be done in a similar way as in \cite{Bernardin/Patricia}. The details are left to the reader. 
	\end{remark}
	
	To prove Theorem \ref{teorema:1}, we first prove a compactness result in Section \ref{sec:4}, more precisely, we prove the tightness of the sequence $(\mathbb{Q}_{N,\theta})_{N\geq 1}$. We then prove in Section \ref{sec:5} that all the limit points of $(\mathbb{Q}_{N,\theta})_{N \geq 1}$ are concentrated on measures that are absolutely continuous with respect to the Lebesgue measure, and whose Radon-Nikodym derivatives are weak solutions to \eqref{eq:1.6}. To such an end we will require some auxiliary results that are given in Sections \ref{sec:6} and \ref{sec:7}. Finally, uniqueness of the weak solutions are obtained in Appendix \ref{appendix}, which in turn gives us the weak convergence and convergence in probability of $(\mathbb{Q}_{N,\theta})_{N \geq 1}$ to $\mathbb{Q}^*_\theta$ as $ N \to \infty $, where the measure $\mathbb{Q}^*_\theta$ is a Dirac measure concentrated on $(\rho,\varrho)du$, and $(\rho,\varrho)$ is the unique weak solution to \eqref{eq:1.6}.
	
	\section{Tightness}\label{sec:4}
	In this section, we show that the sequence of probability measures $(\mathbb{Q}_{N,\theta})_N$ is tight in the Skorohod space $\mathcal{D}([0,T],\, \mathcal{M}_+\times \mathcal{M}^d)$. In order to do that, we invoke the Aldous' criterion and from \cite[Chapter 4, Proposition 1.7]{Kipnis/Landim} it is enough to show that for every function $G$ in a dense subset of $C(D^d)$ with respect to the uniform topology, the sequence of measures that correspond to the real processes $\langle \pi^{k,N}_{t,\theta}, G\rangle$, is tight for each $k \in \{0,1,\dots,d\}$. To such an end, we need to prove two things. First:
	\begin{equation}\label{eq:1.27}
		\lim_{A \to +\infty}\lim_{N \to +\infty}\mathbb{P}_{\mu^N}\Big( |\langle \pi^{k,N}_{t,\theta},G\rangle|>A\Big)=0, \mbox{ for each } k \in \{0,1,\dots ,d\}.
	\end{equation}
	But this follows from Chebychev's inequality and the fact that for the exclusion type dynamics, the number of particles per site is at most one for each fixed velocity. Second, we need to prove that for all $\varepsilon > 0$ and any function $G$ in a dense subset of $C(D^d)$, with respect to the uniform topology:
	\begin{equation}\label{eq:1.28}
		\lim_{\delta \to 0}\limsup_{N \to +\infty}\,\sup_{\stackrel{\tau \in \mathfrak{T}_T}{t \leq \delta}}\mathbb{P}_{\mu^N}\Big(\eta:  |\langle \pi^{k,N}_{\tau+t,\theta},G\rangle-\langle \pi^{k,N}_{t,\theta},G\rangle|>\varepsilon\Big)=0,
	\end{equation}
	where $\mathfrak{T}_T$ denotes the set of stopping times with respect to the canonical filtration bounded by $T$.
	Recall that it is enough to prove the assertion for functions $G$ in a dense subset of $C(D^d)$ with respect to the uniform topology. We will split the proof into two cases, the first one for $\theta \geq 1$ and the second one for $0 \leq \theta < 1 $.
	
	By Dynkin's formula, see, for example, \cite[Appendix A1, Lemma 5.1]{Kipnis/Landim}, for each $k \in \{0,1,\dots,d\}$
	\begin{equation}\label{eq:1.13}
		M^{N,k}_{t,\theta}(G)= \langle \pi^{k,N}_{t,\theta}, G \rangle-\langle 		\pi^{k,N}_{0,\theta}, G \rangle-\int_0^t (\mathcal{L}_{N,\theta}+\partial_s)\langle \pi^{k,N}_{s,\theta}, G \rangle\,ds
	\end{equation}  
	is a martingale with respect to the natural filtration $\mathcal{F}_{t,\theta} = \sigma(\eta_{s,\theta}, s \leq t)$. Then, given $\tau\in\mathfrak{T}_T$, 
	\begin{equation*} \begin{split}
			&\mathbb{P}_{\mu^N} \Big(\eta:  \Big|\langle \pi^{k,N}_{\tau+t,\theta},G \rangle-\langle \pi^{k,N}_{\tau,\theta},G\rangle|>\varepsilon\Big) \\
			=\,& \mathbb{P}_{\mu^N} \Big(\eta:  | M^{N,k}_{\tau,\theta}(G)-M^{N,k}_{\tau+ t,\theta}(G)+\displaystyle\int_{\tau}^{\tau+ t}\mathcal{L}_{N,\theta}\langle \pi^{k,N}_{s,\theta}, G\rangle \,ds\Big|>\varepsilon\Big)\\
			\leq\, & \mathbb{P}_{\mu^N} \Big(\eta:  | M^{N,k}_{\tau,\theta}(G)-M^{N,k}_{\tau+t,\theta}(G)\Big|>\displaystyle\frac{\varepsilon}{2}\Big)+ \mathbb{P}_{\mu^N} \Big(\eta: \Big|\displaystyle\int_{\tau}^{\tau+t}\mathcal{L}_{N,\theta}\langle \pi^{k,N}_{s,\theta}, G\rangle \,ds\Big|>\frac{\varepsilon}{2}\Big).
	\end{split}\end{equation*}
	Applying Chebychev's inequality (resp. Markov's inequality) in the first (resp. second) term  on the right-hand side of the last inequality, we can bound the previous expression from above by
	\begin{equation}\nonumber
		\displaystyle\frac{4}{{\varepsilon}^2}\mathbb{E}_{\mu^N} \left[\Big(M^{N,k}_{\tau,\theta}(G)-M^{N,k}_{\tau+t,\theta}(G)\Big)^2\right]+\frac{2}{\varepsilon}\mathbb{E}_{\mu^N} \left[\, \Big|\int_{\tau}^{\tau+ t}\mathcal{L}_{N,\theta}\langle \pi^{k,N}_{s,\theta}, G\rangle \, ds\Big|\, \right].
	\end{equation}
	Therefore, in order to prove \eqref{eq:1.28} it is enough to show that 
	\begin{equation}\label{eq:1.29}
		\displaystyle\lim_{\delta \to 0} \limsup_{N \to +\infty}\, \sup_{\stackrel{\tau \in \mathfrak{T}_T}{t\leq \delta}} \mathbb{E}_{\mu^N}\Big[\,\Big|\int_{\tau}^{\tau+t}\mathcal{L}_{N,\theta}\langle \pi^{k,N}_{s,\theta}, G\rangle\, ds\Big|\,\Big] = 0, \mbox{~and~}
	\end{equation}
	\begin{equation}\label{eq:1.30}
		\displaystyle\lim_{\delta \to 0} \limsup_{N \to +\infty}\, \sup_{\stackrel{\tau \in \mathfrak{T}_T}{t \leq \delta}} \mathbb{E}_{\mu^N} \Big[\Big( M^{N,k}_{\tau,\theta}(G) - M^{N,k}_{\tau + t,\theta}(G) \Big)^2\Big] = 0.
	\end{equation}
	Let us start by proving \eqref{eq:1.29} for functions $G \in C^2(D^d)$, since by a standard $L^1$ procedure it is easy to extend it to functions $G \in C(D^d)$. Now, we will show that there exists a constant $C$ such that $\mathcal{L}_{N,\theta}\langle \pi^{k,N}_{s,\theta},G\rangle \leq C$ for any $s \leq T$. For that purpose, note that 
	\begin{equation*}
		|\mathcal{L}_{N,\theta}\langle \pi^{k,N}_{s,\theta},G \rangle | \leq |N^2\mathcal{L}^{ex,1}_N\langle \pi^{k,N}_{s,\theta},G\rangle|+|N^2\mathcal{L}^{ex,2}_N\langle \pi^{k,N}_{s,\theta},G\rangle|+|N^2\mathcal{L}^{c}_N\langle \pi^{k,N}_{s,\theta},G\rangle|+|N^2\mathcal{L}^{b}_{N,\theta}\langle \pi^{k,N}_{s,\theta},G\rangle|.
	\end{equation*}
	Let us bound this separately. Note that, 
	\begin{equation}
		\begin{split}\label{eq:1.31}
			|N^2\mathcal{L}^{ex,1}_N\langle \pi^{k,N}_{s,\theta},G\rangle| & \leq \Big|\langle \pi^{k,N}_{s,\theta},\displaystyle\frac{1}{2}\Delta_N G\rangle\Big|+\Big| \displaystyle\frac{N}{2N^d}\sum_{\stackrel{x \in D^d_N}{x_1=1}}\sum_{v \in \mathcal{V}}v_k\eta_{s,\theta}(1,\tilde{x},v)\partial_{x_1}^{N,+}G\left(\tfrac{0}{N},\tfrac{\tilde{x}}{N}\right)\Big|\\
			& +\Big| \displaystyle\frac{N}{2N^d}\sum_{\stackrel{x \in D^d_N}{x_1=N-1}}\sum_{v \in \mathcal{V}}v_k\eta_{s,\theta}(N-1,\tilde{x},v)\partial_{x_1}^{N,-}G\left(\tfrac{N}{N},\tfrac{\tilde{x}}{N}\right)\Big|\\
			&\leq \displaystyle\frac{C}{2}\|G''\|_{\infty} +\frac{C NN^{d-1}}{2N^d}\|G'\|_{\infty}+\frac{C NN^{d-1}}{2N^d}\|G'\|_{\infty} = \displaystyle\frac{C}{2}\|G''\|_{\infty}+C\|G'\|_{\infty},
		\end{split}
	\end{equation}
	where \begin{equation*}\begin{split}
			&\partial_{x_1}^{N,+}G\left(\tfrac{0}{N},\tfrac{\tilde{x}}{N}\right)=N\left[G\left(\tfrac{1}{N},\tfrac{\tilde{x}}{N}\right)-G\left(\tfrac{0}{N},\tfrac{\tilde{x}}{N}\right)\right],\\&\partial_{x_1}^{N,-}G\left(\tfrac{N}{N},\tfrac{\tilde{x}}{N}\right)=-N\left[G\left(\tfrac{N-1}{N}, \tfrac{\tilde{x}}{N}\right)-G\left(\tfrac{N}{N},\tfrac{\tilde{x}}{N}\right)\right].\end{split}
	\end{equation*}
	
	Similarly, using the fact that $p$ has finite range, we have that
	\begin{equation}
		\begin{split}\label{eq:1.33}
			|N^2\mathcal{L}^{ex,2}_N\langle \pi^{k,N}_{s,\theta},G\rangle| 
			& \leq \Big|\displaystyle\frac{1}{N^d} \sum_{j=1}^{d}\sum_{x \in D^d_N}(\partial_{x_j}^NG)\left(\tfrac{x}{N}\right)\sum_{v \in \mathcal{V}}v_k\sum_{z \in \mathbb{Z}^d}p(z,v)z_j\eta_{s,\theta}(x,v)(1-\eta_{s,\theta}(x+z,v))\Big|\\
			&\leq \displaystyle\frac{\tilde{C}N^d\|G'\|_{\infty}}{N^d}
			= \tilde{C}\|G'\|_{\infty}.
	\end{split}\end{equation}
	Also,
	\begin{equation}
		\begin{split}\label{eq:1.32}
			|N^2\mathcal{L}^{b}_{N,\theta}\langle \pi^{k,N}_{s,\theta},G\rangle|&\leq \Big| \displaystyle\frac{N^2}{2N^dN^{\theta}}\sum_{\stackrel{x \in D^d_N}{x_1=1}}\sum_{v \in \mathcal{V}}v_kG\left(\tfrac{1}{N},\tfrac{\tilde{x}}{N}\right)[\hat{\alpha}_v\left(\tfrac{\tilde{x}}{N}\right)-\eta_{s,\theta}(1,\tilde{x},v)]\Big|\\
			& +\Big| \displaystyle\frac{N^2}{2N^dN^{\theta}}\sum_{\stackrel{x \in D^d_N}{x_1=N-1}}\sum_{v \in \mathcal{V}}v_kG\left(\tfrac{N-1}{N},\tfrac{\tilde{x}}{N}\right)[\hat{\beta}_v\left(\tfrac{\tilde{x}}{N}\right)-\eta_{s,\theta}(N-1,\tilde{x},v)]\Big|\\
			& \leq \displaystyle\frac{CN^2N^{d-1}}{2N^dN^{\theta}}\|G\|_{\infty}+ \displaystyle\frac{CN^2N^{d-1}}{N^dN^{\theta}}\|G\|_{\infty}=CN^{1-\theta}\|G\|_{\infty}.
		\end{split}
	\end{equation}
	Now observe that,
	\begin{equation*}
		N^2\mathcal{L}_N^{c}\langle \pi^{k,N}_{s,\theta}, G \rangle = \displaystyle\frac{N^2}{N^d}\sum_{x \in D^d_N}\sum_{v \in \mathcal{V}}v_k G\left(\displaystyle\tfrac{x}{N}\right)\mathcal{L}_N^{c}(\eta_{s,\theta}(x,v)).
	\end{equation*}
	Since the operator is linear, we just need to compute $\mathcal{L}_N^{c}(\eta_{s,\theta}(x,v))$. For $f(\eta)=\eta(x,v)$, we have that 
	\begin{equation*}
		(\mathcal{L}_N^{c}f)(\eta) = \displaystyle\sum_{y \in D^d_N}\sum_{q \in Q}p_c(y,q,\eta)[f(\eta^{y,q})-f(\eta)] = \displaystyle\sum_{q \in Q}p_c(x,q,\eta)[\eta(x,v_{i+2})-\eta(x,v)]=0.
	\end{equation*}
	Therefore,
	\begin{equation}\label{eq.c}
		N^2\mathcal{L}_N^{c}\langle \pi^{k,N}_{s,\theta}, G \rangle =0.
	\end{equation}
	Now, we restrict to the case $\theta \geq 1$. From \eqref{eq:1.31}, \eqref{eq:1.33}, \eqref{eq:1.32} and \eqref{eq.c} we have that $$|\mathcal{L}_{N,\theta}\langle \pi^{k,N}_{s,\theta},G\rangle | \leq C$$ which proves \eqref{eq:1.29}. 
	
	In the case, $0\leq \theta < 1$, if we try to apply the same strategy used for $\theta \geq 1$ we will run into trouble when trying to control the modulus of continuity of $\int_0^t N^2\mathcal{L}^{b}_{N,\theta}\langle \pi^{k,N}_{s,\theta},G\rangle \, ds $, because the expression in \eqref{eq:1.32} can explode when $N \to \infty$. We will prove \eqref{eq:1.28} first for functions $G \in C^2_c(D^d)$ and then we can extend it, by a $L^1$ approximation procedure which is explained below, to functions $G \in C^1(D^d)$. In this case it holds
	\begin{equation}\label{eq:ex1c2}
		\begin{split}
			|N^2\mathcal{L}^{ex,1}_N\langle \pi^{k,N}_{s,\theta},G\rangle| & \displaystyle\leq |\langle \pi^{k,N}_{s,\theta},\frac{1}{2}\Delta_N G\rangle|+\Big| \displaystyle\frac{N}{2N^d}\sum_{\stackrel{x \in D^d_N}{x_1=1}}\sum_{v \in \mathcal{V}}v_k\eta_{s,\theta}(1,\tilde{x},v)\partial_{x_1}^{N,+}G(0,\tilde{x})\Big|\\
			& +\Big| \displaystyle\frac{N}{2N^d}\sum_{\stackrel{x \in D^d_N}{x_1=N-1}}\sum_{v \in \mathcal{V}}v_k\eta_{s,\theta}(N-1,\tilde{x},v)\partial_{x_1}^{N,-}G(1,\tilde{x})\Big|\leq \displaystyle\frac{1}{2}\|G''\|_{\infty}.
		\end{split}
	\end{equation}
	Also,
	\begin{equation}\begin{split}\label{eq:bound}
			|N^2\mathcal{L}^b_N\langle \pi^{k,N}_{s,\theta},G\rangle | 
			&\leq \Big| \displaystyle\frac{N^2}{2N^dN^{\theta}}\sum_{\stackrel{x \in D^d_N}{x_1=1}}\sum_{v \in \mathcal{V}}v_kG\left(\tfrac{1}{N},\tfrac{\tilde{x}}{N}\right)[\hat{\alpha}_v\left(\tfrac{\tilde{x}}{N}\right)-\eta_{s,\theta}(1,\tilde{x},v)]\Big|\\
			& +\Big| \displaystyle\frac{N^2}{2N^dN^{\theta}}\sum_{\stackrel{x \in D^d_N}{x_1=N-1}}\sum_{v \in \mathcal{V}}v_kG\left(\tfrac{N-1}{N},\tfrac{\tilde{x}}{N}\right)[\hat{\beta}_v\left(\tfrac{\tilde{x}}{N}\right)-\eta_{s,\theta}(N-1,\tilde{x},v)]\Big|\\
			&=0,
	\end{split}\end{equation}
	for $N$ sufficiently large. This is a consequence of the fact that the function $G$ has compact support in the interior of $D^d$. Using \eqref{eq:1.33}, \eqref{eq.c}, \eqref{eq:ex1c2} and \eqref{eq:bound} this finishes the proof of \eqref{eq:1.29} for any $\theta \geq 0$. 
	
	We will now prove \eqref{eq:1.30}. Since $$\Big( M^{N,k}_{\tau,\theta}(G) - M^{N,k}_{\tau + t,\theta}(G) \Big)^2 - \displaystyle\int_{\tau}^{\tau+t} \mathcal{L}_{N,\theta}[\langle \pi^{k,N}_{s,\theta},G\rangle]^2-2\langle \pi^{k,N}_{s,\theta},G\rangle  \mathcal{L}_{N,\theta}\langle \pi^{k,N}_{s,\theta},G\rangle ds$$ is a mean zero martingale, we obtain that
	$$
	\mathbb{E}_{\mu^N} \Big[\Big( M^{N,k}_{\tau,\theta}(G) - M^{N,k}_{\tau + t,\theta}(G) \Big)^2\Big]= \mathbb{E}_{\mu^N} \Big[\int_{\tau}^{\tau+t} \mathcal{L}_{N,\theta}[\langle \pi^{k,N}_{s,\theta},G\rangle]^2-2\langle \pi^{k,N}_{s,\theta},G\rangle  \mathcal{L}_{N,\theta}\langle \pi^{k,N}_{s,\theta},G\rangle ds\Big].
	$$  
	Note that, \eqref{eq:1.30} holds if we show that
	$$
	\displaystyle\int_{\tau}^{\tau+t} \mathcal{L}_{N,\theta}[\langle \pi^{k,N}_{s,\theta},G\rangle]^2-2\langle \pi^{k,N}_{s,\theta},G\rangle  \mathcal{L}_{N,\theta}\langle \pi^{k,N}_{s,\theta},G\rangle ds,
	$$
	converges to zero uniformly in $t \in [0,T]$, when $N \to \infty$.
	Simple, albeit long, computations show that
	\begin{equation*}\begin{split}
			&N^2\mathcal{L}_N^{ex,1}\langle \pi^{k,N}_{s,\theta},G\rangle^2-2\langle \pi^{k,N}_{s,\theta},G\rangle N^2 \mathcal{L}_N^{ex,1}\langle \pi^{k,N}_{s,\theta},G\rangle\\
			& = \displaystyle\frac{1}{2N^{2d}}\sum_{v \in \mathcal{V}}\sum_{x \in D^d_N}\sum_{j=1}^dv_k^2\, \big[\eta_{s,\theta}(x,v)-\eta_{s,\theta}(x+e_j,v)\big]^2[\partial^N_{x_j}G\left(\tfrac{x}{N}\right)]^2\leq \displaystyle\frac{C}{N^d}\|G'\|^2_{\infty}.
	\end{split}\end{equation*}
	We also have
	\begin{equation*}
		\begin{split}
			&N^2\mathcal{L}_N^{ex,2}\langle \pi^{k,N}_{s,\theta},G\rangle^2-2\langle \pi^{k,N}_{s,\theta},G\rangle N^2 \mathcal{L}_N^{ex,2}\langle \pi^{k,N}_{s,\theta},G\rangle\\
			& = \displaystyle\frac{1}{N^{2d+1}}\sum_{v \in \mathcal{V}}\sum_{x \in D^d_N}\sum_{w \in \mathbb{Z}^d}v_k^2\,\eta_{s,\theta}(x,v)(1-\eta_{s,\theta}(x+w,v))p(w,v)w_j^2 [\partial_{x_j}^NG\left(\tfrac{x}{N}\right)]^2\leq \displaystyle\frac{\tilde{C}}{N^{d+1}}\|G'\|^2_{\infty}.
		\end{split}
	\end{equation*}
	Moreover,   
	\begin{equation*}
		\begin{split}
			&	N^2\mathcal{L}_{N,\theta}^{b}\langle \pi^{k,N}_{s,\theta},G\rangle^2-2\langle \pi^{k,N}_{s,\theta},G\rangle N^2 \mathcal{L}_{N,\theta}^{b}\langle \pi^{k,N}_{s,\theta},G\rangle\\
			&	= \displaystyle\frac{N^2}{2N^{2d}}\displaystyle\sum_{\stackrel{x \in D^d_N}{x_1=1}}\sum_{v \in \mathcal{V}}\left[\frac{\hat{\alpha}_v(\tfrac{\tilde{x}}{N})(1-\eta_{s,\theta}(x,v))+(1-\hat{\alpha}_v(\tfrac{\tilde{x}}{N}))\eta_{s,\theta}(x,v)}{N^{\theta}}\right]v_k^2\,[G\left(\tfrac{1}{N},\tfrac{\tilde{x}}{N}\right)]^2 \\
			&+ \displaystyle\frac{N^2}{2N^{2d}}\displaystyle\sum_{\stackrel{x \in D^d_N}{x_1=N-1}}\sum_{v \in \mathcal{V}}\left[\frac{\hat{\beta}_v(\tfrac{\tilde{x}}{N})(1-\eta_{s,\theta}(x,v))+(1-\hat{\beta}_v(\tfrac{\tilde{x}}{N}))\eta_{s,\theta}(x,v)}{N^{\theta}}\right]v_k^2\,[G\left(\tfrac{N-1}{N},\tfrac{\tilde{x}}{N}\right)]^2\leq \displaystyle\frac{\tilde{C}}{N^{d-1}}\|G\|^2_{\infty}.
		\end{split}
	\end{equation*}
	And doing a long, but simple, computation, we can show  that 
	\begin{equation}\nonumber
		N^2\mathcal{L}_N^{c}\langle \pi^{k,N}_{s,\theta},G\rangle^2-2\langle \pi^{k,N}_{s,\theta},G\rangle N^2 \mathcal{L}_N^{c}\langle \pi^{k,N}_{s,\theta},G\rangle=0.
	\end{equation}
	This finishes the proof of tightness for any $\theta \geq 0$.
	
	\section{Characterization of  limit points}\label{sec:5}
	This section deals with the characterization of the limit of points of $(\mathbb{Q}_{N,\theta})_N$ for $\theta \geq 0$. For simplicity we split the proof in three different regimes.
	\subsection{Characterization of the limit points for $\theta \in [0,1)$}\label{subsec:5.1}
	\begin{proposition}\label{prop:<1}
		If $\mathbb{Q}_{\theta}^*$ is a limit point of $(\mathbb{Q}_{N,\theta})_{N\geq 1}$, then
		\begin{equation*}
			\mathbb{Q}_{\theta}^*\left[\pi_{\cdot}: \pi=(\rho,\varrho)du \mbox{ and } \mathfrak{F}_{+\infty}((\rho,\varrho),t,G)=0\right]=1,
		\end{equation*}
		for all $t \in [0,T], \,\forall G \in C^{1,2}_0([0,T]\times D^d)$ and $\theta \in [0,1)$, where the functional $\mathfrak{F}_{\mathfrak{C}}$ was defined in \eqref{eq:1.7}.
	\end{proposition}
	
	\begin{proof}
		It is enough to verify that, for $\delta >0$ and $G \in  C^{1,2}_0([0,T]\times D^d)$ fixed,
		\begin{equation}\label{eq:char_limit_theta01}
			\mathbb{Q}_\theta^*\Bigg[\pi_{\cdot}: \pi^k=\varrho^k\,du \,\mbox{ and } \displaystyle\sup_{0 \leq t \leq T}|\mathfrak{F}^k_{+\infty}((\rho,\varrho),t,G)| > \delta \Bigg]=0,
		\end{equation}
		where $\mathfrak{F}^k_{\mathfrak{C}}$ corresponds to the $k$-th evolution equation in  $\mathfrak{F}_{\mathfrak{C}}$ defined in \eqref{eq:1.7} for $k \in \{0,1,\dots,d\}$.
		First, observe that 
		$$\mathbb{Q}_\theta^*\Big[\pi_{\cdot}: \pi^k=\varrho^k\,du\Big] = 1,$$
		for $k\in \{0,1,\dots,d\}$. Now, we need to replace the term $\chi(\theta_v(\Lambda(\rho,\varrho)))$ in \eqref{eq:1.7} by a functional of the empirical measure. To this end, for each $\varepsilon>0$, let $j_\varepsilon(\cdot)$, which we will use to obtain an approximation of the point valuation, be given by
		$$
		j_\varepsilon(x) = (2\varepsilon)^{-d} \mathbbm{1}_{[-\varepsilon,\varepsilon]^d}(x),\quad x\in D^d.
		$$
		
		Now, let $(\rho,\varrho)\ast j_\varepsilon (\cdot) = \left(\rho\ast j_\varepsilon (\cdot), \varrho_1\ast j_\varepsilon (\cdot), \ldots, \varrho_d\ast j_\varepsilon (\cdot)\right)$. We have that $(\rho,\varrho)\ast j_\varepsilon (\cdot) \to (\rho,\varrho)(\cdot)$ in $L^1(D^d)$ as $\varepsilon \to 0$. Since $\chi(\theta_v(\Lambda(\cdot)))$ is Lipschitz, we have that $\chi(\theta_v(\Lambda((\rho,\varrho)\ast j_\varepsilon (\cdot)))) \to \chi(\theta_v(\Lambda((\rho,\varrho)(\cdot))))$ in $L^1(D^d)$. Thus, $\mathbb{Q}^\ast_{\theta}$ almost surely, for every $v\in\mathcal{V}$,
		$$
		\displaystyle\int_0^t\int_{D^d}\tilde{v}\, \chi(\theta_v(\Lambda((\rho, \varrho)\ast j_\varepsilon)))\sum_{i=1}^{d}v_i\frac{\partial G}{\partial u_i}(r,u) \,du\,dr \to \displaystyle\int_0^t\int_{D^d}\tilde{v}\, \chi(\theta_v(\Lambda(\rho, \varrho)))\sum_{i=1}^{d}v_i\frac{\partial G}{\partial u_i}(r,u) \,du\,dr,
		$$
		as $\varepsilon \to 0$. Let, now, $\mathfrak{F}_{\mathfrak{C}}^\varepsilon((\rho,\varrho),t,G)$ be the expression obtained after we replace $\chi(\theta_v(\Lambda(\rho,\varrho)))$ by $\chi(\theta_v(\Lambda((\rho,\varrho)\ast j_\varepsilon)))$ in \eqref{eq:1.7}. Therefore, \eqref{eq:char_limit_theta01} will follow if we show that
		\begin{equation}\label{eq:char_limit_theta01_with_epsilon}
			\limsup_{\varepsilon\to 0}\mathbb{Q}_\theta^*\Bigg[\pi_{\cdot}: \pi^k=\varrho^k\,du \,\mbox{ and } \displaystyle\sup_{0 \leq t \leq T}|\mathfrak{F}^{\varepsilon,k}_{+\infty}((\rho,\varrho),t,G)| > \delta \Bigg]=0,
		\end{equation}	
		where $\mathfrak{F}^{\varepsilon,k}_{\mathfrak{C}}$ corresponds to the $k$-th projection of $\mathfrak{F}^{\varepsilon}_{\mathfrak{C}}$ defined above, for each $k \in \{0,1,\dots,d\}$.
		By the same arguments as in \cite[p.77]{Kipnis/Landim}, we have that
		\begin{equation*}
			\textbf{I}^{[N\varepsilon]}(x,\eta) = C_{N,\varepsilon} (\pi^N_{\theta}\ast j_\varepsilon)\left(\tfrac{x}{N}\right),
		\end{equation*}
		where $\pi^N_{\theta}\ast j_\varepsilon(\cdot) = \left(\pi_{\theta}^{0,N}\ast j_\varepsilon(\cdot),\pi_{\theta}^{1,N}\ast j_\varepsilon(\cdot),\ldots,\pi_{\theta}^{d,N}\ast j_\varepsilon(\cdot) \right)$ and $C_{N,\varepsilon} = 1 + O(N^{-1})$.
		Since the set considered in \eqref{eq:char_limit_theta01_with_epsilon} is an open set with respect to Skorohod topology in $\mathcal{D}([0,T],\mathcal{M}_+ \times \mathcal{M}^d)$, we can use the Portmanteau's Theorem directly and bound the probability in \eqref{eq:char_limit_theta01_with_epsilon} by
		\begin{equation}\label{eq:carac}
			\begin{split}
				&\limsup_{\varepsilon\to 0}\displaystyle\liminf_{N \to \infty} \mathbb{Q}_{N,\theta}\left[\pi_{\cdot}: \displaystyle\sup_{0 \leq t \leq T}\left| \langle \pi^{k,N}_t , G_t \rangle -\langle \pi^{k,N}_0 , G_0 \rangle \right. \right.-\displaystyle\int_0^t \left\langle \pi^{k,N}_r, \left(\partial_r G(r,u)+\frac{1}{2}\Delta G_r\right)\right\rangle \,dr\\
				&+\displaystyle\frac{1}{2}\int_0^t\int_{\{1\}\times\mathbb{T}^{d-1}}\sum_{v \in \mathcal{V}}v_k\beta_v(\tilde{u})\frac{\partial G}{\partial u_1}(r,1,\tilde{u})\,dS\,dr-\frac{1}{2}\int_0^t\int_{\{0\}\times \mathbb{T}^{d-1}} \sum_{v \in \mathcal{V}}v_k\alpha_v(\tilde{u})\frac{\partial G}{\partial u_1}(r,0,\tilde{u})\,dS\,dr \\
				&-\left.\left.\displaystyle\int_0^t\int_{D^d}\sum_{v \in \mathcal{V}}v_k\, \chi(\theta_v(\Lambda(\textbf{I}^{[N\varepsilon]}(u,\eta))))\sum_{i=1}^{d}v_i\frac{\partial G}{\partial u_i}(r,u) \,du\,dr\right| > \delta \right],
			\end{split}
		\end{equation}
		for each $k \in \{0,1,\dots,d\}$. Summing and subtracting $\displaystyle\int_0^t\mathcal{L}_{N,\theta}\langle  \pi^{k,N}_r, G \rangle \,dr$ to the expression inside the supremum in \eqref{eq:carac}, we can bound it by the sum of
		\begin{equation*}
			\displaystyle\liminf_{N \to \infty} \mathbb{P}_{\mu^N}\left[\eta_\cdot: \sup_{0 \leq t \leq T}|M^{N,k}_t(G)| > \frac{\delta}{2}\right] \mbox{~and~}
		\end{equation*}
		\begin{equation}\label{eq.conti}
			\begin{split}
				\displaystyle\limsup_{\varepsilon\to 0}\limsup_{N \to \infty}\mathbb{P}_{\mu^N}&\left[\eta_{\cdot}:\displaystyle \sup_{0 \leq t \leq T}\left|\int_0^t \mathcal{L}_{N,\theta}\langle \pi^{k,N}_r , G_r \rangle \, dr - \frac{1}{2}\displaystyle\int_0^t \langle \pi^{k,N}_r, \Delta G_r \rangle \,dr\right.\right.\\ 
				&-\displaystyle\int_0^t\int_{D^d}\sum_{v \in \mathcal{V}}v_k\, \chi(\theta_v(\Lambda(\textbf{I}^{[N\varepsilon]}(u,\eta))))\sum_{i=1}^{d}v_i\frac{\partial G}{\partial u_i}(r,u) \,du\,dr\\
				&+\displaystyle\frac{1}{2}\int_0^t\int_{\{1\}\times\mathbb{T}^{d-1}}\frac{\partial G}{\partial u_1}(r,1,\tilde{u})\sum_{v \in \mathcal{V}}v_k \hat\beta_v(\tilde{u})\,dS\,dr\\
				&-\left.\left.\frac{1}{2}\int_0^t\int_{\{0\}\times \mathbb{T}^{d-1}}\frac{\partial G}{\partial u_1}(r,0,\tilde{u}) \sum_{v \in \mathcal{V}}v_k \hat\alpha_v(\tilde{u})\,dS\,dr\right| > \displaystyle\frac{\delta}{2}\right].
			\end{split}
		\end{equation}
		By a long computation we have that
		\begin{equation}\label{eq:intger}
			\begin{split}
				\displaystyle\int_0^t \mathcal{L}_{N,\theta} \langle \pi^{k,N}_r, G_r\rangle \, dr & =\displaystyle\frac{1}{2}\int_0^t\langle \pi^{k,N}_r, \Delta_NG_r \rangle \, dr\\
				& - \displaystyle\frac{1}{2} \int_0^t \frac{1}{N^{d-1}} \sum_{\stackrel{x \in D^d_N}{x_1=N-1}} I_k(\eta_r(N-1,\tilde{x}))\partial^{N,-}_{x_1}G_r\left(\tfrac{N}{N},\tfrac{\tilde{x}}{N}\right)\, dr\\ 
				& + \displaystyle\frac{1}{2}\int_0^t\frac{1}{N^{d-1}} \sum_{\stackrel{x \in D^d_N}{x_1=1}} I_k(\eta_r(1,\tilde{x})) \partial^{N,+}_{x_1}G_r\left(\tfrac{0}{N},\tfrac{\tilde{x}}{N}\right)\, dr\\ 
				&+\displaystyle\int_0^t \displaystyle\frac{1}{N^{d}}\displaystyle\sum_{j=1}^d \sum_{x \in D^d_N}(\partial^N_{x_j}G_r)\left(\tfrac{x}{N}\right)\tau_x W^{N,\eta_r}_{j,k}\, dr\\ 
				& +\displaystyle\int_0^t \displaystyle\frac{N^2}{2N^{d+\theta}}\sum_{\stackrel{x \in D^d_N}{x_1=1}}  G_r\left(\tfrac{1}{N},\tfrac{\tilde{x}}{N}\right) \Big[\sum_{v \in \mathcal{V}}v_k\hat{\alpha}_v(\tfrac{\tilde{x}}{N})- I_k(\eta_r(1,\tilde{x}))\Big]\,dr\\ 
				& +\displaystyle\int_0^t\displaystyle\frac{N^2}{2N^{d+\theta}}\sum_{\stackrel{x \in D^d_N}{x_1=N-1}} G_r\left(\tfrac{N-1}{N}, \tfrac{\tilde{x}}{N}\right) \Big[\sum_{v \in \mathcal{V}}v_k \hat{\beta}_v(\tfrac{\tilde{x}}{N})-I_k(\eta_r(N-1,\tilde{x}))\Big]\,dr,
			\end{split}
		\end{equation}
		where $\tau_x$ stands for the translation by $x\in D^d_N$, so that $(\tau_x\eta)(y,v)=\eta(x+y,v)\mathbbm{1}_{\{x+y \in D^d_N\}}$ for all $y \in \mathbb{Z}^d, v \in \mathcal{V},$ and $W^{N,\eta_r}_{j,k}$ is given by:
		$$
		W^{N,\eta_r}_{j,k}=\displaystyle\sum_{v \in \mathcal{V}}v_k\sum_{z \in \mathbb{Z}^d}p(z,v)z_j\eta_r(0,v)[1-\eta_r(z,v)],
		$$
		where we recall that $v_0=1$ and $\pi^{k,N}_t$ is the empirical measure defined in \eqref{eq:1.11}.
		Now, restricting \eqref{eq:intger} for the case $G \in C^{1,2}_0([0,T]\times D^d)$, we obtain that
		\begin{equation}\label{eq:intgerDBC}
			\begin{split}
				\displaystyle\int_0^t \mathcal{L}_{N,\theta} \langle \pi^{k,N}_r, G_r\rangle \, dr & =\displaystyle\frac{1}{2}\int_0^t\langle \pi^{k,N}_r, \Delta_NG_r \rangle \, dr\\
				& - \displaystyle\frac{1}{2} \int_0^t \frac{1}{N^{d-1}} \sum_{\stackrel{x \in D^d_N}{x_1=N-1}} I_k(\eta_r(N-1,\tilde{x}))\partial^{N,-}_{x_1}G_r\left(\tfrac{N}{N},\tfrac{\tilde{x}}{N}\right)\, dr\\ 
				& + \displaystyle\frac{1}{2}\int_0^t\frac{1}{N^{d-1}} \sum_{\stackrel{x \in D^d_N}{x_1=1}} I_k(\eta_r(1,\tilde{x})) \partial^{N,+}_{x_1}G_r\left(\tfrac{0}{N},\tfrac{\tilde{x}}{N}\right)\, dr\\ 
				&+\displaystyle\int_0^t \displaystyle\frac{1}{N^{d}}\displaystyle\sum_{j=1}^d \sum_{x \in D^d_N}(\partial^N_{x_j}G_r)\left(\tfrac{x}{N}\right)\tau_x W^{N,\eta_r}_{j,k}\, dr
			\end{split}
		\end{equation}
		
		Using the expression \eqref{eq:intgerDBC} in \eqref{eq.conti} let us bound the expression inside the probability in equation \eqref{eq.conti} by the sum of the following terms
		\begin{subequations}\label{eqs}
			\begin{align}
				\displaystyle\sup_{0 \leq t\leq T} &  \Big|  \int_0^t \Big[ \frac{1}{2N^d}\sum_{x \in D^d_N}I_k(\eta_r(x))\Delta_NG_r\left(\tfrac{x}{N}\right) - \frac{1}{2N^d}\sum_{x \in D^d_N}I_k(\eta_r(x))\Delta G_r\left(\tfrac{x}{N}\right) \Big] \, dr\Big|,\label{eq1} \\
				\displaystyle\sup_{0 \leq t \leq T}  & \Big| \int_0^t \frac{1}{2N^{d-1}}\sum_{\stackrel{x \in D^d_N}{x_1=N-1}} I_k(\eta_r(N-1,\tilde{x}))\Big[\partial_{x_1}^{N,-}G_r\left(\tfrac{N}{N},\tfrac{\tilde{x}}{N}\right)-\partial_{x_1}G_r(\tfrac {1}{N}, \tfrac{ \tilde{x}}{N})\Big]\, dr\Big|, \label{eq2} \\ 
				\displaystyle\sup_{0 \leq t \leq T} & \Big| \frac{1}{2} \int_0^t \frac{1}{N^{d-1}}\sum_{\stackrel{x \in D^d_N}{x_1=N-1}}\partial_{u_1}G_r(\tfrac{1}{N}, \tfrac{\tilde{x}}{N})\Big[ I_k(\eta_r(N-1,\tilde{x})) -\sum_{v \in \mathcal{V}}v_k \beta_v(\tfrac{\tilde{x}}{N}) \Big]\, dr\Big|, \label{eq3} \\
				\displaystyle\sup_{0 \leq t \leq T}  & \Big| \int_0^t \frac{1}{2N^{d-1}}\sum_{\stackrel{x \in D^d_N}{x_1=1}} I_k(\eta_r(1,\tilde{x}))\Big[\partial_{x_1}^{N,+}G_r\left(\tfrac{0}{N},\tfrac{\tilde{x}}{N}\right)-\partial_{x_1}G_r(\tfrac{0}{N},\tfrac{\tilde{x}}{N})\Big]\, dr\Big|, \label{eq4}\\
				\displaystyle\sup_{0 \leq t \leq T} & \Big| \frac{1}{2}\int_0^t \frac{1}{N^{d-1}} \sum_{\stackrel{x \in D^d_N}{x_1=1}}\partial_{u_1}G_r(\tfrac{0}{N}, \tfrac{\tilde{x}}{N})\Big[ I_k(\eta_r(1,\tilde{x})) - \sum_{v \in \mathcal{V}}v_k \alpha_v(\tfrac{\tilde{x}}{N})\Big]\, dr\Big|, \label{eq5}\\
				\displaystyle\sup_{0 \leq t \leq T}  &\Big|\int_0^t \Big[ \frac{1}{N^d} \sum_{j=1}^d \sum_{x \in D^d_N}(\partial_{x_j}^N G_r)(\tfrac{x}{N})\tau_xW_{j,k}-\int_{D^d}\sum_{v \in \mathcal{V}}v_k\, \chi(\theta_v(\Lambda(\textbf{I}^{[N\varepsilon]}(\tfrac{x}{N},\eta))))\sum_{i=1}^d v_i\frac{\partial G_r}{\partial u_i} \, dx \Big] dr\Big| \label{eq8}
			\end{align}
		\end{subequations}
		
		Using  the Taylor expansion, we obtain that \eqref{eq1}, \eqref{eq2} and \eqref{eq4} converges to zero as $N \to \infty$.  From the Replacement Lemmas (resp. Proposition \ref{prop:2} and Lemma \ref{lema:5}) it follows that \eqref{eq3}, \eqref{eq5} and  \eqref{eq8} vanishes as $N\to+\infty$ and then $\varepsilon\to 0$. This concludes the proof.
	\end{proof}
	
	\subsection{Characterization of limit points for $\theta = 1$}\label{subsec:7.2}
	We begin by fixing some notations used along this subsection. For $k\in\{0,1,\cdots, d\}$, we denote by
	\begin{equation}\label{overrighteta}
		\begin{array}{lcl}
			\overrightarrow{\eta}^{k,\ell}(x_1,\tilde{x},v):=\displaystyle\frac{1}{\ell}\sum_{y \in \overrightarrow{\Lambda}_{x_1}^{\ell}}I_k(\eta_s(y,\tilde{x},v)) & \mbox{ and }&	\overleftarrow{\eta}^{k,\ell}(x_1,\tilde{x},v):=\displaystyle\frac{1}{\ell}\sum_{y \in \overleftarrow{\Lambda}_{x_1}^\ell} I_k(\eta_s(y,\tilde{x},v)),
		\end{array}
	\end{equation}
	where $\overleftarrow{\Lambda}_{x_1}^{\ell}:=\{x_1-\ell+1,\dots,x_1\}$ (resp. $\overrightarrow{\Lambda}_{x_1}^{\ell}:=\{x_1,\dots,x_1+\ell-1\}$).
	\begin{proposition}\label{prop:=1}
		If $\mathbb{Q}_\theta^*$ is a limit point of $(\mathbb{Q}_{N,\theta})_{N\geq 1}$, then 
		\begin{equation}\nonumber
			\mathbb{Q}_\theta^*\left[\pi_{\cdot}: \pi=(\rho,\varrho)du \mbox{ and } \mathfrak{F}_{1}((\rho,\varrho),t,G)=0\right]=1,
		\end{equation}
		for all $t \in [0,T], \,\forall G \in C^{1,2}([0,T]\times D^d)$ and $\theta=1$.
	\end{proposition}
	\begin{proof}
		It is enough to verify that, for $\delta >0$ and $G \in C^{1,2}([0,T]\times D^d)$ fixed,
		\begin{equation*}
			\mathbb{Q}_\theta^*\left[\pi_{\cdot}: \pi^k=\varrho^k\,du \mbox{ and } \sup_{0 \leq t \leq T}\left|\mathfrak{F}^k_{1}((\rho,\varrho),t,G)
			\right| > \delta \right]=0,
		\end{equation*}
		for all $k\in\{0,1,\cdots, d\}$.
		Since $\mathbb{Q}_\theta^*\Big[\pi_{\cdot}: \pi^k=\varrho^k\,du\Big] = 1$, following the same strategy as for the case $\theta\in [0,1)$, we obtain that it is enough to prove that
		\begin{equation}\label{eq:emevidencia}
			\begin{split}
				&\limsup_{\varepsilon\to 0} \mathbb{Q}_\theta^*\Bigg[\pi_{\cdot}: \displaystyle\sup_{0 \leq t \leq T}\Bigg|  \displaystyle\int_{D^d} \varrho^k(t,u)G(t,u)\, du-\int_{D^d}\varrho^k(0,u)G(0,u)\,du \\ 
				&-\displaystyle\int_0^t\int_{D^d}\sum_{v \in \mathcal{V}}v_k\, \chi(\theta_v(\Lambda((\rho, \varrho)\ast j_\varepsilon))))\sum_{i=1}^{d}v_i\frac{\partial G}{\partial u_i}(r,u) \,du\,dr\\ 
				&-\displaystyle\int_0^t\int_{D^d}\varrho^k(r,u)\left(\partial_r G(r,u)+\frac{1}{2}\Delta G_r\right)\, du\,dr\\ 
				&-\displaystyle\frac{1}{2}\int_0^t \int_{\{1\}\times \mathbb{T}^{d-1}}
				G(r,1,\tilde{u})\sum_{v \in \mathcal{V}}v_k\beta_v(\tilde{u}) \,dS\,dr
				-\displaystyle\frac{1}{2}\int_0^t \int_{\{0\}\times \mathbb{T}^{d-1}}\!\!\!
				G(r,0,\tilde{u})\sum_{v \in \mathcal{V}}v_k\alpha_v(\tilde{u})\,dS\,dr\\ 
				&				+\displaystyle\frac{1}{2}\int_0^t\int_{\{1\}\times\mathbb{T}^{d-1}}\varrho^k(r,1,\tilde{u})\left[\frac{\partial G}{\partial u_1}(r,1,\tilde{u})+G(r,1,\tilde{u}) \right]\,dS\,dr\\ 
				&-\displaystyle\frac{1}{2}\int_0^t\int_{\{0\}\times \mathbb{T}^{d-1}} \varrho^k(r,0,\tilde{u})\left[\frac{\partial G}{\partial u_1}(r,0,\tilde{u})-G(r,0,\tilde{u})\right]\,dS\,dr \Bigg| > \delta \Bigg]=0.
			\end{split}
		\end{equation}
		At this point, we would like to apply Portmanteau's Theorem in the  probabilities $\mathbb{Q}_{N,\theta}$, as we did in the previous case. But here we face a problem. Unfortunately, the set inside the above probability is not an open set with respect to the Skorohod topology. In order to avoid this problem, we fix $\varepsilon >0$ (which is the same $\varepsilon$ used in \eqref{eq:emevidencia}) and we consider two approximations of the identity, for fixed $u_1 \in [0,1]$ which are given on $w \in [0,1]$ by
		\begin{equation*}
			\overrightarrow{\imath}^{0}_{\varepsilon}(w)=\displaystyle\frac{1}{\varepsilon}\mathbbm{1}_{[0,\varepsilon)}(w) \mbox{~~and~~}\overleftarrow{\imath}^{1}_{\varepsilon}(w)=\displaystyle\frac{1}{\varepsilon}\mathbbm{1}_{(1-\varepsilon,1]}(w).
		\end{equation*}
		We use the notation $\langle \pi^k_r, \overrightarrow{\imath}^{0}_{\varepsilon} \rangle  = \langle \varrho^k_r, \overrightarrow{\imath}^{0}_{\varepsilon} \rangle$ and  $\langle \pi^k_r, \overleftarrow{\imath}^{1}_{\varepsilon} \rangle = \langle \varrho^k_r, \overleftarrow{\imath}^{1}_{\varepsilon} \rangle$, so that
		\begin{equation*}
			\langle \pi^k_r, \overrightarrow{\imath}^{0}_{\varepsilon} \rangle(\tilde{u}) = \displaystyle\frac{1}{\varepsilon} \int_{0}^{\varepsilon} \varrho^k_r(w, \tilde{u})\,dw \mbox{~~and~~} \langle \pi^k_r, \overleftarrow{\imath}^{1}_{\varepsilon} \rangle (\tilde{u}) =  \displaystyle\frac{1}{\varepsilon} \int_{1-\varepsilon}^{1} \varrho^k_r(w, \tilde{u})\,dw.
		\end{equation*}
		By summing and subtracting proper terms, we bound the probability in \eqref{eq:emevidencia} from above by 
		\begin{equation}\label{eq:emevidencia2}
			\begin{split}
				&\limsup_{\varepsilon\to 0}\mathbb{Q}_\theta^*\left[\pi_{\cdot}: \displaystyle\sup_{0 \leq t \leq T}\left|  \displaystyle\int_{D^d} \varrho^k(t,u)G(t,u)\, du-\int_{D^d}\varrho^k(0,u)G(0,u)\,du \right.\right.\\ 
				&-\displaystyle\int_0^t\int_{D^d}\sum_{v \in \mathcal{V}}v_k \chi(\theta_v(\Lambda((\rho, \varrho)\ast j_\varepsilon)))\sum_{i=1}^{d}v_i\frac{\partial G}{\partial u_i}(r,u) \,du\,dr\\ 
				&+\displaystyle\frac{1}{2}\int_0^t\int_{\{1\}\times\mathbb{T}^{d-1}}\left[\varrho^k(r,1,\tilde{u})-\langle \pi^k_r, \overleftarrow{\imath}^{1}_{\varepsilon}\rangle(\tilde{u}) \right]\left[\frac{\partial G}{\partial u_1}(r,1,\tilde{u})+G(r,1,\tilde{u}) \right]\,dS\,dr\\ 
				&+\displaystyle\frac{1}{2}\int_0^t\int_{\{1\}\times\mathbb{T}^{d-1}}\langle \pi^k_r, \overleftarrow{\imath}^{1}_{\varepsilon}\rangle(\tilde{u})\left[\frac{\partial G}{\partial u_1}(r,1,\tilde{u})+G(r,1,\tilde{u}) \right]\,dS\,dr\\ 
				&-\displaystyle\frac{1}{2}\int_0^t\int_{\{0\}\times \mathbb{T}^{d-1}} \left[\varrho^k(r,0,\tilde{u})  - \langle \pi^k_r, \overrightarrow{\imath}^{0}_{\varepsilon} \rangle(\tilde{u}) \right]\left[\frac{\partial G}{\partial u_1}(r,0,\tilde{u})-G(r,0,\tilde{u})\right]\,dS\,dr \\ 
				&-\displaystyle\frac{1}{2}\int_0^t\int_{\{0\}\times \mathbb{T}^{d-1}} \langle \pi^k_r, \overrightarrow{\imath}^{0}_{\varepsilon} \rangle (\tilde{u})\left[\frac{\partial G}{\partial u_1}(r,0,\tilde{u})-G(r,0,\tilde{u})\right]\,dS\,dr \\ 
				&-\displaystyle\frac{1}{2}\int_0^t \int_{\{1\}\times \mathbb{T}^{d-1}}
				G(r,1,\tilde{u})\sum_{v \in \mathcal{V}}v_k\beta_v(\tilde{u})  \,dS\,dr	-\displaystyle\frac{1}{2}\int_0^t \int_{\{0\}\times \mathbb{T}^{d-1}}
				G(r,0,\tilde{u})\sum_{v \in \mathcal{V}}v_k\alpha_v(\tilde{u}) \,dS\,dr\\ 
				&-\left.\left.\displaystyle\int_0^t\int_{D^d}\varrho^k(r,u)\left(\partial_r G(r,u)+\frac{1}{2}\Delta G_r\right)\, du\,dr\right| > \delta \right]=0.
			\end{split}	
		\end{equation}
		From Lebesgue's Differentiation Theorem, observe that 
		\begin{equation*}
			\lim_{\varepsilon \to 0} \left|\varrho^k(r,0,\tilde{u}) - \langle \pi^k_r, \overrightarrow{\imath}^{0}_{\varepsilon} \rangle(\tilde{u}) \right|=0
			\mbox{~~and~~}	\lim_{\varepsilon \to 0} \left|\varrho^k(r,1,\tilde{u})-\langle \pi^k_r, \overleftarrow{\imath}^{1}_{\varepsilon}\rangle(\tilde{u}) \right|=0.
		\end{equation*}
		Since the functions $\overrightarrow{\imath}^{0}_{\varepsilon}$ and $\overleftarrow{\imath}^{1}_{\varepsilon}$ are not continuous, we cannot use Portmanteau's Theorem. However, we can approximate each one of these functions by continuous functions, in such a way that the error vanishes as $\varepsilon \to 0$. Then, since the set inside the probability in \eqref{eq:emevidencia2} is an open set with respect to the Skorohod topology, we can use Portmanteau's Theorem and bound it from above by
		\begin{equation}\label{eq:emevidencia3}
			\begin{split}
				&\displaystyle\limsup_{\varepsilon\to 0}\liminf_{N \to \infty}\mathbb{Q}_{N,\theta}\Bigg[\pi_{\cdot}: \displaystyle\sup_{0 \leq t \leq T}\Bigg|   \langle \pi^{k}_t,G_t \rangle -\langle \pi^k_0, G_0\rangle -\displaystyle\int_0^t \langle \pi^{k}_r,(\partial_r +\frac{1}{2}\Delta) G_r\rangle\,dr\\ 
				&-\displaystyle\int_0^t\int_{D^d}\sum_{v \in \mathcal{V}}v_k\, \chi(\theta_v(\Lambda(\textbf{I}^{[N\varepsilon]}(u,\eta))))\sum_{i=1}^{d}v_i\frac{\partial G}{\partial u_i}(r,u) \,du\,dr\\ 
				&-\displaystyle\frac{1}{2}\int_0^t \int_{\{0\}\times \mathbb{T}^{d-1}} G(r,1,\tilde{u})\sum_{v \in \mathcal{V}}v_k\beta_v(\tilde{u})  \,dS\,dr -\displaystyle\frac{1}{2}\int_0^t \int_{\{0\}\times \mathbb{T}^{d-1}}
				G(r,0,\tilde{u})\sum_{v \in \mathcal{V}}v_k\alpha_v(\tilde{u}) \,dS\,dr\\ 
				&+\displaystyle\frac{1}{2}\int_0^t \sum_{\stackrel{u \in D^d_N}{u_1=N-1}}\overleftarrow{\eta}^{k,\varepsilon N}_r(N-1,\tilde{u},v)\left[\frac{\partial G}{\partial u_1}(r,1,\tilde{u})+G(r,1,\tilde{u})\right]\,dS\,dr \\ 
				&-\displaystyle\frac{1}{2}\int_0^t \sum_{\stackrel{u \in D^d_N}{u_1=1}}\overrightarrow{\eta}^{k,\varepsilon N}_r(1,\tilde{u},v) \left[\frac{\partial G}{\partial u_1}(r,0,\tilde{u})-G(r,0,\tilde{u})\right]\,dS\,dr \Bigg| > \frac{\delta}{2} \Bigg]=0. 
			\end{split}
		\end{equation}
		Summing and subtracting $\displaystyle\int_0^t \mathcal{L}_{N,\theta}\langle \pi^{k,N}_r,G_r\rangle\, dr$ in equation inside the supremum in \eqref{eq:emevidencia3}, we can bound it from above by the sum of
		\begin{equation}\label{eq.anula}
			\mathbb{P}_{\mu^N}\left[\sup_{0 \leq t \leq T}|M^{N,k}_t(G)| > \frac{\delta}{4}\right] \mbox{~and~}
		\end{equation} 
		\begin{equation}	\label{eq:main}	
			\begin{split}
				&\mathbb{P}_{\mu^N}\left[\displaystyle\sup_{0 \leq t \leq T}\left| 	\displaystyle\int_0^t \mathcal{L}_{N,\theta}\langle \pi^{k,N}_r,G_r\rangle\, dr-\displaystyle\int_0^t\int_{D^d}\sum_{v \in \mathcal{V}}v_k\, \chi(\theta_v(\Lambda(\textbf{I}^{[N\varepsilon]}(u,\eta))))\sum_{i=1}^{d}v_i\frac{\partial G}{\partial u_i}(r,u) \,du\,dr\right.\right.\\ 
				&+\displaystyle\frac{1}{2}\int_0^t\sum_{\stackrel{u \in D^d_N}{u_1=N-1}} \overleftarrow{\eta}^{k,\varepsilon N}_r(N-1,\tilde{u},v)\left[\frac{\partial G}{\partial u_1}(r,1,\tilde{u})+G(r,1,\tilde{u})\right]\,dS\,dr \\ 
				&-\displaystyle\frac{1}{2}\int_0^t \sum_{\stackrel{u \in D^d_N}{u_1=1}}\overrightarrow{\eta}^{k,\varepsilon N}_r(1,\tilde{u},v) \left[\frac{\partial G}{\partial u_1}(r,0,\tilde{u})-G(r,0,\tilde{u})\right]\,dS\,dr \\ 
				&-\displaystyle\frac{1}{2}\int_0^t \int_{\{1\}\times \mathbb{T}^{d-1}} G(r,1,\tilde{u})\sum_{v \in \mathcal{V}}v_k\beta_v(\tilde{u})  \,dS\,dr
				-\displaystyle\frac{1}{2}\int_0^t \int_{\{0\}\times \mathbb{T}^{d-1}}
				G(r,0,\tilde{u})\sum_{v \in \mathcal{V}}v_k\alpha_v(\tilde{u}) \,dS\,dr\\ 
				&-\left.\left.\displaystyle\int_0^t \langle \pi^k_r,\frac{1}{2}\Delta G_r\rangle\,dr\right| > \frac{\delta}{4} \right]=0.
			\end{split}
		\end{equation}
		From Doob's inequality, the equation \eqref{eq.anula} vanishes as $N \to \infty$. 
		On the other hand using the value of $\int_0^t \mathcal{L}_{N,\theta}\langle \pi^{k,N}_r,G_r\rangle\, dr$ (see \eqref{eq:intger}), using the same argument from the previous subsection we can bound last display from above by a sum of \eqref{eq1}, \eqref{eq8}, in addition to this terms that we have already obtained previously, we will have
		\begin{equation}\label{eq:final}\begin{split}
				&\sup_{0 \leq t \leq T} \left| \int_0^t \frac{1}{2N^{d-1}}\left[ I_k(\eta_r(1,\tilde{u})) - \overrightarrow{\eta}_r^{k,\varepsilon N }(1,\tilde{u},v)\right] \tilde{G}(r,\tilde{u}) \,dr\right| \mbox{~~and~~}\\
				&\sup_{0 \leq t \leq T} \left| \int_0^t \frac{1}{2N^{d-1}}\left[ I_k(\eta_r(N-1,\tilde{u}))- \overleftarrow{\eta}_r^{k,\varepsilon N}(N-1,\tilde{u},v)\right] \tilde{G}(r,\tilde{u}) \,dr\right|.
		\end{split}\end{equation}
		Using the Replacement Lemmas it is easy to see that equations \eqref{eq:final} converges to zero, as $N \to \infty$ and $\varepsilon \to 0$. This concludes the proof.
	\end{proof}

	\subsection{Characterization of limit points for $\theta > 1$}\label{subsec:7.3}
	\begin{proposition}\label{prop:>1}
		If $\mathbb{Q}_{\theta}^*$ is a limit point of $\{\mathbb{Q}_{N,\theta}\}_{N\geq 1}$, then 
		\begin{equation}\nonumber
			\mathbb{Q}_{\theta}^*\left[\pi_{\cdot}: \pi=(\rho,\varrho)du \mbox{ and } \mathfrak{F}_{0}((\rho,\varrho),t,G)=0\right]=1,
		\end{equation}
		for all $t \in [0,T], \,\forall G \in C^{1,2}([0,T]\times D^d)$ and $\theta>1$.
	\end{proposition}
	\begin{proof}
		The proof follows the same reasoning used in the proof of the previous Propositions \ref{prop:<1} and \ref{prop:=1} and it will be omitted here.
	\end{proof}
	\section{Replacement Lemmas}\label{sec:6}
	We will need the following auxiliary function $h$ to be able to obtain some entropy estimates that are essential to the proof of the hydrodynamic limit. 
	For each $v \in \mathcal{V}$, consider the functions $h^v_k:D^d \to (0,1)$, for $k\in\{0,\dots,d\}$. 
	\begin{remark}\label{obs:3}
		We will have two situations for the function $h= \displaystyle\sum_{v \in \mathcal{V}}(h^v_0,v_1h^v_1,\dots , v_d h^v_d)$:
		\begin{itemize}
			\item [(1)] when $\theta \in [0,1)$ we will assume that for each $k \in \{0,\dots,d\}$, $h^v_k$ are smooth functions and that there exists $\delta > 0$, such that the restriction of $h$ to $[0,\delta)\times \mathbb{T}^{d-1}$ equals to the vector valued function $\alpha(\tilde{u})$, and that the restriction of $h$ to $(1-\delta,1]\times \mathbb{T}^{d-1}$ equals to the vector valued function $\beta(\tilde{u})$.
			\item [(2)] when $\theta \geq 1$ we will assume that $h$ is a constant function.
		\end{itemize}
	\end{remark}
	
	We then consider $\nu^N_{h}$ as the product measure on $X_N$ with marginals given by
	\begin{equation}\label{eq:1.12}
		\nu^N_{h}\{\eta: \eta(x, \cdot) = \xi\} = m_{\Lambda(h(x/N))}(\xi),
	\end{equation}
	where $m_{\lambda}(\cdot)$ was defined in \eqref{eq:1.2}.

	We observe that for  $h$ satisfying (1) above, for $N$ sufficiently large, we have that 
	\begin{equation}\label{eq:medidas}
		\begin{split}
			&\nu_h^N(\eta:\eta(1,\tilde x,v)=1)=m_{\Lambda(h(x/N))}(1)=\hat{\alpha}_v(\tilde x/N)\\
			&\nu_h^N(\eta:\eta(N-1,\tilde x,v)=1)=m_{\Lambda(h(x/N))}(0)=\hat{\beta}_v(\tilde x/N).
		\end{split}
	\end{equation}

	\subsection{Estimates on Dirichlet Forms}\label{subsec:6.1}
	Let $f: X_N \to \mathbb{R}$ be a local function. 
	Note that $\langle \mathcal{L}_{N,\theta}\sqrt{f},\sqrt{f}\rangle_{\nu^N_{h}}$ does not always have a closed form. To estimate $\langle \mathcal{L}_{N,\theta}\sqrt{f},\sqrt{f}\rangle_{\nu^N_{h}}$,  let
	\begin{equation*}
		D_{\nu^N_{h}}^\theta(\sqrt{f}) = D^{ex}_{\nu^N_{h}}(\sqrt{f}) + D^{c}_{\nu^N_{h}}(\sqrt{f}) + D^{b,\theta}_{\nu^N_{h}}(\sqrt{f}), 
	\end{equation*}
	with
	\begin{equation*}
		D^{ex}_{\nu^N_{h}}(\sqrt{f}) = \displaystyle\sum_{v \in \mathcal{V}}\sum_{x \in D^d_N}\sum_{z \in D^d_N}c_{(x,z,v)}(\eta)\int[\sqrt{f(\eta^{x,z,v})}-\sqrt{f(\eta)}]^2 d\nu^N_{h},
	\end{equation*}
	where $c_{(x,z,v)}(\eta)=\eta(x,v)(1-\eta(z,v)) P_N(z-x,v)$,
	\begin{equation*}
		D^{c}_{\nu^N_{h}}(\sqrt{f}) = \displaystyle\sum_{q \in \mathcal{Q}}\sum_{x \in D^d_N}\int p_c(x,q,\eta)[\sqrt{f(\eta^{x,q})}-\sqrt{f(\eta)}]^2 d\nu^N_{h}
	\end{equation*}
	where $\mathcal{Q} \mbox{ and } p_c(x,q,\eta)$ were defined in \eqref{def:conjQ} and \eqref{eq:defpc}, respectively, and
	\begin{equation*}
		\begin{split}
			D^{b,\theta}_{\nu^N_{h}}(\sqrt{f}) & = \displaystyle\sum_{v \in \mathcal{V}}\sum_{x \in \{1\}\times \mathbb{T}^{d-1}_N}\int \frac{r^N_{x}(\eta,\alpha)}{2N^{\theta}}[\sqrt{f(\sigma^{x,v}\eta)}-\sqrt{f(\eta)}]^2 d\nu^N_{h} \\
			& + \displaystyle\sum_{v \in \mathcal{V}}\sum_{x \in \{N-1\}\times \mathbb{T}^{d-1}_N}\int \frac{r^N_{x}(\eta,\beta)}{2N^{\theta}} [\sqrt{f(\sigma^{x,v}\eta)}-\sqrt{f(\eta)}]^2 d\nu^N_{h} := D^{b,\theta,\alpha}_{\nu^N_{h}}(\sqrt{f}) + D^{b,\theta,\beta}_{\nu^N_{h}}(\sqrt{f}),
		\end{split}
	\end{equation*}
	where $r^N_{x}(\eta,\alpha)$ and $r^N_{x}(\eta,\beta)$ were defined in \eqref{eq:taxas}. In order to prove the main result of this section, we need some intermediate results. For that purpose, we adapt from \cite{Bernardin/Patricia} the following lemmas:
	\begin{lemma}\label{lema:1}
		Let $T: \eta \in X_N \to T(\eta) \in X_N$ be a transformation in the configuration space and $c:\eta \in X_N \to c(\eta)$ be a positive local function.
		Let $f$ be a density with respect to a probability measure $\nu^N_{h}$ on $X_N$. Then, 
		\begin{equation*}
			\begin{split}
				&	\langle c(\eta)[\sqrt{f(T(\eta))}-\sqrt{f(\eta)}], \sqrt{f(\eta)}\rangle_{\nu^N_{h}} \leq
				-\displaystyle\frac{1}{4}\int c(\eta)\left(\sqrt{f(T(\eta))}-\sqrt{f(\eta)}\right)^2\,d\nu^N_{h}\\
				&+\displaystyle\frac{1}{16}\int \frac{1}{c(\eta)}\left[c(\eta)-c(T(\eta))\frac{\nu^N_{h}(T(\eta))}{\nu^N_{h}(\eta)}\right]^2\left(\sqrt{f(T(\eta))}+\sqrt{f(\eta)}\right)^2 \,d\nu^N_{h}.
			\end{split}
		\end{equation*}
	\end{lemma}
	
	\begin{lemma}\label{lema:2}
		Let $f$ be a density with respect to a probability measure $\nu^N_{h}$ on $X_N$. Then, 
		\begin{equation}\nonumber
			\begin{array}{l}
				\displaystyle\sup_{x\neq y}\int f(\eta^{x,y,v}) d\nu^N_{h} \leq C, \quad
				\displaystyle\sup_{x}\int f(\eta^{x,q}) d\nu^N_{h} \leq C \quad\textrm{and}\quad
				\displaystyle\sup_{x}\int f(\sigma^{x,v}\eta) d\nu^N_{h} \leq C,
			\end{array}
		\end{equation}
		where $\eta^{x,y,v}$, $\eta^{x,q}$ and $\sigma^{x,v}\eta$, were defined in \eqref{eq:config}, \eqref{eq:configc} and \eqref{eq:configb}, respectively.
	\end{lemma}
	\begin{remark}
		Note that in \cite{Bernardin/Patricia} the proof is done for the cases $T(\eta)=\eta^{x,y,v}$ and $T(\eta)=\sigma^{x,v}\eta$ in Lemma \ref{lema:2}. For $T(\eta)=\eta^{x,q}$, the proof follows from the arguments of the other cases.
	\end{remark}
	As a consequence of Lemmas \ref{lema:1} and \ref{lema:2}, we conclude that
	\begin{corollary}\label{coro:1}
		Let $h$ be one of the functions satisfying the assumption of Remark \ref{obs:3}. Let $f:X_N \to \mathbb{R}$ be a density with respect to the measure $\nu^N_{h}$. Then,
		\begin{enumerate}
			\item if $h$ is a constant function, 
			$
			N^2\langle \mathcal{L}^{ex}_N\sqrt{f}, \sqrt{f}\rangle_{\nu^N_{h}} = -\displaystyle\frac{N^2}{2}D_{\nu^N_{h}}^{ex}(\sqrt{f});
			$
			\item if $h$ is a smooth function, 
			$
			N^2\langle \mathcal{L}^{ex}_N\sqrt{f}, \sqrt{f}\rangle_{\nu^N_{h}} = -\displaystyle\frac{N^2}{4}D_{\nu^N_{h}}^{ex}(\sqrt{f})+ \mathcal{R}_N(h)
			$
			with $|\mathcal{R}_N(h)|\leq CN^d$.
		\end{enumerate}
	\end{corollary}
	
	\begin{proof}
		Begin by observing that
		\begin{equation*}
			\begin{split}
				N^2\langle \mathcal{L}^{ex}_N\sqrt{f}, \sqrt{f}\rangle_{\nu^N_{h}} & =\displaystyle{N^2}\int\sum_{v \in \mathcal{V}}\sum_{x \in D^d_N}\sum_{z \in D^d_N}c_{(x,z,v)}(\eta)\left(\sqrt{f(\eta^{x,z,v})}-\sqrt{f(\eta)}\right) \sqrt{f(\eta)}\, d{\nu^N_{h}}\\
				& +\displaystyle\frac{N^2}{2}\int\sum_{v \in \mathcal{V}}\sum_{x \in D^d_N}\sum_{z \in D^d_N}c_{(x,z,v)}(\eta)\left(\sqrt{f(\eta^{x,z,v})}-\sqrt{f(\eta)}\right) \sqrt{f(\eta^{x,z,v})}\, d{\nu^N_{h}}\\
				& -\displaystyle\frac{N^2}{2}\int\sum_{v \in \mathcal{V}}\sum_{x \in D^d_N}\sum_{z \in D^d_N}c_{(x,z,v)}(\eta)\left(\sqrt{f(\eta^{x,z,v})}-\sqrt{f(\eta)}\right) \sqrt{f(\eta^{x,z,v})}\, d{\nu^N_{h}},
			\end{split}
		\end{equation*}
		where $c_{(x,z,v)}(\eta)=\eta(x,v)(1-\eta(z,v))P_N(z-x,v)$.
		Combining some terms and doing a change of variables, we obtain that the last display equals to 
		\begin{equation*}
			\begin{split}
				&-\displaystyle\frac{N^2}{2}\int \sum_{v \in \mathcal{V}}\sum_{x \in D^d_N}\sum_{z \in D^d_N}c_{(x,z,v)}(\eta)\left(\sqrt{f(\eta^{x,z,v})}-\sqrt{f(\eta)}\right)^2\, d{\nu^N_{h}}\\
				&+\displaystyle\frac{N^2}{2}\int\sum_{v \in \mathcal{V}}\sum_{x \in D^d_N}\sum_{z \in D^d_N}c_{(x,z,v)}(\eta)\left(\sqrt{f(\eta^{x,z,v})}-\sqrt{f(\eta)}\right)\sqrt{f(\eta^{x,z,v})}\, d{\nu^N_{h}}\\
				&-\displaystyle\frac{N^2}{2}\int\sum_{v \in \mathcal{V}}\sum_{x \in D^d_N}\sum_{z \in D^d_N}c_{(z,x,v)}(\eta)\displaystyle\frac{\nu^N_{h}(\eta^{x,z,v})}{\nu^N_{h}(\eta)}\left(\sqrt{f(\eta^{x,z,v})}-\sqrt{f(\eta)}\right)\sqrt{f(\eta^{x,z,v})} d{\nu^N_{h}}.
			\end{split}	
		\end{equation*}
		Last display equals to
		\begin{equation*}
			\begin{split}
				&-\displaystyle\frac{N^2}{2}\int\sum_{v \in \mathcal{V}}\sum_{x \in D^d_N}\sum_{z \in D^d_N}c_{(x,z,v)}(\eta)\left(\sqrt{f(\eta^{x,z,v})}-\sqrt{f(\eta)}\right)^2\, d{\nu^N_{h}}\\
				&+\displaystyle\frac{N^2}{2}\int\sum_{v \in \mathcal{V}}\sum_{x \in D^d_N}\sum_{z \in D^d_N}c_{(x,z,v)}(\eta)\left(1- \frac{\nu^N_{h}(\eta^{x,z,v})}{\nu^N_{h}(\eta)}\right)\sqrt{f(\eta^{x,z,v})}\left(\sqrt{f(\eta^{x,z,v})}-\sqrt{f(\eta)}\right)\, d{\nu^N_{h}}.
			\end{split}
		\end{equation*}
		Hence, $N^2\langle \mathcal{L}^{ex}_N\sqrt{f}, \sqrt{f}\rangle_{\nu^N_{h}} =-\displaystyle\frac{N^2}{2}D^{ex}_{\nu^N_{h}}(\sqrt{f}) + g_N(h)$, where 
		\begin{equation}\nonumber
			g_N(h)=\! \displaystyle\frac{N^2}{2}\!\int\sum_{v \in \mathcal{V}}\sum_{x \in D^d_N}\!\sum_{z \in D^d_N}\!c_{(x,z,v)}(\eta)\!\! \left(1-\tfrac{\nu^N_{h}(\eta^{x,z,v})}{\nu^N_{h}(\eta)}\right)\!\!\sqrt{f(\eta^{x,z,v})}\left(\sqrt{f(\eta^{x,z,v})}-\sqrt{f(\eta)}\right) d{\nu^N_{h}}.
		\end{equation}
		To handle $g_N(h)$, we start by observing that $\left|1-\displaystyle\frac{\nu^N_{h}(\eta^{x,z,v})}{\nu^N_{h}(\eta)}\right| = \left|1-\displaystyle\frac{\gamma_{z,v}}{\gamma_{x,v}}\right|$ where
		\begin{equation}\label{eq:1.16}
			\gamma_{x,v}=\displaystyle\frac{\theta_v(\Lambda(h(x/N)))}{1-\theta_v(\Lambda(h(x/N)))}.
		\end{equation}
		Thus, if $h$ is constant, then $g_N(h)=0$. 
		
		On the other hand, let us assume that $h$ is not constant, we need to redo the  analysis of $g_N(h)$. By applying the elementary inequality, $ab \leq \tfrac{1}{2}a^2+\tfrac{1}{2}b^2$ in $g_N(h)$, with 
		$$
		a=\frac{N}{\sqrt{2}}\sqrt{\eta(x,v)(1-\eta(z,v)P_N(z-x,v))}\left(\sqrt{f(\eta^{x,z,v})}-\sqrt{f(\eta)}\right)
		$$ 
		and 
		$$
		b=\frac{N}{\sqrt{2}}\sqrt{\eta(x,v)(1-\eta(z,v)P_N(z-x,v))}\sqrt{f(\eta^{x,z,v})}\left(1-\frac{\nu^N_{h}(\eta^{x,z,v})}{\nu^N_{h}(\eta)}\right),
		$$
		we can bound $g_N(h)$ from above by
		\begin{equation*}
			\begin{split}
				& \displaystyle\frac{N^2}{4}\int\sum_{v \in \mathcal{V}}\sum_{x \in D^d_N}\sum_{z \in D^d_N}\eta(x,v)(1-\eta(z,v))P_N(z-x,v)\left(\sqrt{f(\eta^{x,z,v})}-\sqrt{f(\eta)}\right)^2\, d{\nu^N_{h}} \\
				+ &\displaystyle\frac{N^2}{4}\int\sum_{v \in \mathcal{V}}\sum_{x \in D^d_N}\sum_{z \in D^d_N}\eta(x,v)(1-\eta(z,v))P_N(z-x,v) \left(1-\frac{\nu^N_{h}(\eta^{x,z,v})}{\nu^N_{h}(\eta)}\right)^2 f(\eta^{x,z,v})\, d{\nu^N_{h}}.
			\end{split}
		\end{equation*}
		Thus $|g_N(h)|\leq \displaystyle\frac{N^2}{4}D^{ex}_{\nu^N_{h}}(\sqrt{f})+\mathcal{R}_N(h)$, where
		\begin{equation*}
			\mathcal{R}_N(h):= \displaystyle\frac{N^2}{4}\int\sum_{v \in \mathcal{V}}\sum_{x \in D^d_N}\sum_{z \in D^d_N}c_{(x,z,v)}(\eta)\left(1-\frac{\nu^N_{h}(\eta^{x,z,v})}{\nu^N_{h}(\eta)}\right)^2 f(\eta^{x,z,v})\, d{\nu^N_{h}}.
		\end{equation*}
		Doing again the change of variables $\eta^{x,z,v}=\xi$, we obtain 
		\begin{equation*}
			\mathcal{R}_N(h)= \displaystyle\frac{N^2}{4}\int\sum_{v \in \mathcal{V}}\sum_{x \in D^d_N}\sum_{z \in D^d_N}c_{(z,x,v)}(\eta) \left|1-\frac{\nu^N_{h}(\eta)}{\nu^N_{h}(\eta^{x,z,v})}\right|^2  \frac{\nu^N_{h}(\eta^{x,z,v})}{\nu^N_{h}(\eta)}f(\eta)\, d{\nu^N_{h}}.
		\end{equation*}
		Now, observe that
		\begin{equation}\label{eq:1.18}
			\left|1-\frac{\nu^N_{h}(\eta)}{\nu^N_{h}(\eta^{x,z,v})}\right|= \left(1-\frac{\gamma_{x,v}}{\gamma_{z,v}}\right)  \leq \tilde{c}\|\gamma '\|_{\infty}\frac{1}{N},
		\end{equation}
		since $\gamma$ is bounded away from zero, see \eqref{eq:1.16}, and
		\begin{equation}\nonumber
			\left|\frac{\nu^N_{h}(\eta^{x,z,v})}{\nu^N_{h}(\eta)}\right| \leq C.
		\end{equation}
		Also note that $f$ is a density with respect to $\nu^N_{h}$, therefore,
		\begin{equation*}
			| \mathcal{R}_N(h) | \leq CN^d.
		\end{equation*}
		This finishes the proof of Corollary \ref{coro:1}.
	\end{proof}
	\begin{corollary}\label{coro:2}
		For $h$ being  one of the functions satisfying the assumption of Remark \ref{obs:3}  and for a density $f$, with respect to the measure $\nu^N_{h}$, it holds
		\begin{equation}\label{eq:1.19}
			N^2\langle \mathcal{L}^c_N\sqrt{f}, \sqrt{f}\rangle_{\nu^N_{h}} = -\displaystyle\frac{N^2}{2}D_{\nu^N_{h}}^c(\sqrt{f}).
		\end{equation}
	\end{corollary}
	\begin{proof}
		Let $q=(v,w,v',w')$ and $\tilde{q}=(v',w',v,w)$, recall \eqref{eq:defpc}. We have that
		\begin{equation*}
			\begin{split}
				N^2\langle \mathcal{L}^c_N\sqrt{f}, \sqrt{f}\rangle_{\nu^N_{h}} & =   \displaystyle N^2 \int\sum_{y \in D^d_N}\sum_{q \in Q}p_c(y,q,\eta)\left(\sqrt{f(\eta^{y,q})}-\sqrt{f(\eta)}\right)\sqrt{f(\eta)}\,d\nu^N_{h}\\
				& + \displaystyle\frac{N^2}{2}\displaystyle\int\sum_{y \in D^d_N}\sum_{q \in Q}p_c(y,q,\eta)\left(\sqrt{f(\eta^{y,q})}-\sqrt{f(\eta)}\right)\sqrt{f(\eta^{y,q})}\,d\nu^N_{h}\\
				& - \displaystyle\frac{N^2}{2}\displaystyle\int\sum_{y \in D^d_N}\sum_{q \in Q}p_c(y,q,\eta)\left(\sqrt{f(\eta^{y,q})}-\sqrt{f(\eta)}\right)\sqrt{f(\eta^{y,q})}\,d\nu^N_{h}.
			\end{split}
		\end{equation*}
		From a change of variables 
		\begin{equation*}\begin{split}
				N^2\langle \mathcal{L}^c_N\sqrt{f}, \sqrt{f}\rangle_{\nu^N_{h}}  = &  \displaystyle-\frac{N^2}{2} \displaystyle\int\sum_{y \in D^d_N}\sum_{q \in Q}p_c(y,q,\eta)\left(\sqrt{f(\eta^{y,q})}-\sqrt{f(\eta)}\right)^2\,d\nu^N_{h}\\
				& + \displaystyle\frac{N^2}{2}\displaystyle\int\sum_{y \in D^d_N}\sum_{q \in Q}p_c(y,q,\eta)\left(\sqrt{f(\eta^{y,q})}-\sqrt{f(\eta)} \right)\sqrt{f(\eta^{y,q})}\,d\nu^N_{h}\\
				& - \displaystyle\frac{N^2}{2}\displaystyle\int\sum_{y \in D^d_N}\sum_{q \in Q}p_c(y,q,\eta)\left(\sqrt{f(\eta^{y,q})}-\sqrt{f(\eta)} \right)\,\sqrt{f(\eta^{y,q})}\, \frac{\nu^N_{h}(\eta^{y,q})}{\nu^N_{h}(\eta)}\,d\nu^N_{h}.
		\end{split}	\end{equation*}
		Combining some terms, we obtain 
		\begin{equation*}
			\begin{split}
				N^2\langle \mathcal{L}^c_N\sqrt{f}, \sqrt{f}\rangle_{\nu^N_{h}} &= 
				-\displaystyle\frac{N^2}{2} \displaystyle\int\sum_{y \in D^d_N}\sum_{q \in Q}p_c(y,q,\eta)\left(\sqrt{f(\eta^{y,q})}-\sqrt{f(\eta)}\right)^2\,d\nu^N_{h}\\
				& + \displaystyle\frac{N^2}{2}\displaystyle\int\sum_{y \in D^d_N}\sum_{q \in Q}p_c(y,q,\eta)\left(\sqrt{f(\eta^{y,q})}-\sqrt{f(\eta)} \right)\sqrt{f(\eta^{y,q})}\left[1- \frac{\nu^N_{h}(\eta^{y,q})}{\nu^N_{h}(\eta)}\right]\,d\nu^N_{h}.
			\end{split}
		\end{equation*}
		Since $v+w=v'+w'$, we observe that 
		$
		\displaystyle\frac{\nu^N_{h}(\eta^{y,q})}{\nu^N_{h}(\eta)}=\displaystyle\frac{\gamma_{y,v'}\gamma_{y,w'}}{\gamma_{y,v}\gamma_{y,w}}=1,
		$
		therefore,
		$
		N^2\langle \mathcal{L}^c_N\sqrt{f}, \sqrt{f}\rangle_{\nu^N_{h}} = -\displaystyle\frac{N^2}{2}D_{\nu^N_{h}}^c(\sqrt{f}),
		$
		and this identity holds for functions $h$ satisfying the assumption of Remark \ref{obs:3}.
	\end{proof}
	
	\begin{corollary}\label{coro:3}
		For a function $h$ satisfying one of the assumptions of Remark \ref{obs:3} and for a density $f$ with respect to the measure $\nu^N_{h}$ it holds 
		\begin{equation}\label{eq:1.21}
			N^2\langle \mathcal{L}^{b}_{N,\theta}\sqrt{f}, \sqrt{f}\rangle_{\nu^N_{h}} = -\displaystyle\frac{N^2}{2}D_{\nu^N_{h}}^{b,\theta}(\sqrt{f})+ \mathcal{R}^{\vartheta}_{N,\theta}(h) ,
		\end{equation}
		with
		$|\mathcal{R}^{\vartheta}_{N,\theta}(h)|\lesssim \displaystyle\frac{ N^{d+1}}{N^{\theta}}$ if $\theta\geq 1$, while 
		$\mathcal{R}^{\vartheta}_{N,\theta}(h)=0$, for $N$
		sufficiently large, if $\theta\in[0,1)$. Above $\vartheta\in\{\ell, r\}$.
	\end{corollary}
	
	\begin{proof}
		We present the proof for the left boundary since the other case is analogous. Recall \eqref{eq:gen_boundary}. By splitting the integral on the left-hand side of \eqref{eq:1.21} (with respect to the left boundary dynamics) into the integral over the sets $A_0=\{\eta:\eta((1,\tilde{x}),v)=0\}$ and $A_1=\{\eta:\eta((1,\tilde{x}),v)=1\}$, we write it as
		\begin{equation*}\begin{split}
				&\displaystyle\frac{N^2}{2N^{\theta}}\int_{A_0}\sum_{\stackrel{x \in D^d_N}{x_1=1}}\sum_{v \in \mathcal{V}}\hat{\alpha}_v\left(\tfrac{\tilde{x}}{N}\right)[\sqrt{f(\sigma^{x,v}\eta)}-\sqrt{f(\eta)}]\sqrt{f(\eta)}\, d\nu^N_{h}\\
				+&\displaystyle\frac{N^2}{2N^{\theta}}\int_{A_1}\sum_{\stackrel{x \in D^d_N}{x_1=1}}\sum_{v \in \mathcal{V}}\left(1-\hat{\alpha}_v\left(\tfrac{\tilde{x}}{N}\right)\right)[\sqrt{f(\sigma^{x,v}\eta)}-\sqrt{f(\eta)}]\sqrt{f(\eta)}\, d\nu^N_{h}.
		\end{split}	\end{equation*}
		Last display can be rewritten as
		\begin{equation*}
			\begin{split}
				& \displaystyle\frac{N^2}{2N^{\theta}}\int_{A_0}\sum_{\stackrel{x \in D^d_N}{x_1=1}}\sum_{v \in \mathcal{V}}\hat{\alpha}_v\left(\tfrac{\tilde{x}}{N}\right)\left(\sqrt{f(\sigma^{x,v}\eta)}\sqrt{f(\eta)}-\frac{1}{2}\left(\sqrt{f(\eta)}\right)^2\right)\, d\nu^N_{h}\\
				+&\displaystyle\frac{N^2}{2N^{\theta}}\int_{A_1}\sum_{\stackrel{x \in D^d_N}{x_1=1}}\sum_{v \in \mathcal{V}}\left(1-\hat{\alpha}_v\left(\tfrac{\tilde{x}}{N}\right)\right)\left(\sqrt{f(\sigma^{x,v}\eta)}\sqrt{f(\eta)}-\frac{1}{2}\left(\sqrt{f(\eta)}\right)^2\right)\, d\nu^N_{h}\\
				-&\displaystyle\frac{N^2}{4N^{\theta}}\int_{A_0}\sum_{\stackrel{x \in D^d_N}{x_1=1}}\sum_{v \in \mathcal{V}}\hat{\alpha}_v\left(\tfrac{\tilde{x}}{N}\right)\left(\sqrt{f(\eta)}\right)^2\, d\nu^N_{h}\\
				-&\displaystyle\frac{N^2}{4N^{\theta}}\int_{A_1}\sum_{\stackrel{x \in D^d_N}{x_1=1}}\sum_{v \in \mathcal{V}}\left(1-\hat{\alpha}_v\left(\tfrac{\tilde{x}}{N}\right)\right)\left(\sqrt{f(\eta)}\right)^2\, d\nu^N_{h}.
			\end{split}
		\end{equation*}
		Summing and subtracting the term needed to complete the square in last display, we obtain that $N^2\langle \mathcal{L}^{b}_{N,\theta}\sqrt{f}, \sqrt{f}\rangle_{\nu^N_{h}}$ is equal to 
		\begin{equation*}
			\begin{split}
				& -\displaystyle\frac{N^2}{4N^{\theta}}\int_{A_0}\sum_{\stackrel{x \in D^d_N}{x_1=1}}\sum_{v \in \mathcal{V}}\hat{\alpha}_v\left(\tfrac{\tilde{x}}{N}\right)\left(\sqrt{f(\sigma^{x,v}\eta)}-\sqrt{f(\eta)}\right)^2\, d\nu^N_{h}\\
				& -\displaystyle\frac{N^2}{4N^{\theta}}\int_{A_1}\sum_{\stackrel{x \in D^d_N}{x_1=1}}\sum_{v \in \mathcal{V}}\left(1-\hat{\alpha}_v\left(\tfrac{\tilde{x}}{N}\right)\right)\left(\sqrt{f(\sigma^{x,v}\eta)}-\sqrt{f(\eta)}\right)^2\, d\nu^N_{h} \\
				& +\displaystyle\frac{N^2}{4N^{\theta}}\int_{A_0}\sum_{\stackrel{x \in D^d_N}{x_1=1}}\sum_{v \in \mathcal{V}}\hat{\alpha}_v\left(\tfrac{\tilde{x}}{N}\right)\left([\sqrt{f(\sigma^{x,v}\eta)}]^2-[\sqrt{f(\eta)}]^2\right)\, d\nu^N_{h}\\
				& +\displaystyle\frac{N^2}{4N^{\theta}}\int_{A_1}\sum_{\stackrel{x \in D^d_N}{x_1=1}}\sum_{v \in \mathcal{V}}\left(1-\hat{\alpha}_v\right)\left(\tfrac{\tilde{x}}{N}\right)\left([\sqrt{f(\sigma^{x,v}\eta)}]^2-[\sqrt{f(\eta)}]^2\right)\, d\nu^N_{h}. 
			\end{split}
		\end{equation*}
		Using a change of variables on the last two terms above and for a general function $h(\cdot)$, we obtain that the part concerning the left boundary dynamics of $N^2\langle \mathcal{L}^{b}_{N,\theta}\sqrt{f}, \sqrt{f}\rangle_{\nu^N_{h}}$ is equal to
		\begin{equation*}
			-\displaystyle\frac{N^2}{4} D_{\nu^N_{h}}^{b,\theta,\alpha}(\sqrt{f})+	\mathcal{R}^{\alpha}_{N,\theta}(h)
		\end{equation*}
		where 
		\begin{equation*}
			\begin{split}
				\mathcal{R}^{\alpha}_{N,\theta}(h)=&-\displaystyle\frac{N^2}{4N^{\theta}}\int_{A_1}\sum_{\stackrel{x \in D^d_N}{x_1=1}}\sum_{v \in \mathcal{V}}\hat{\alpha}_v\left(\tfrac{\tilde{x}}{N}\right)\left[\frac{1-m_{\Lambda(h(x/N))}(1)}{m_{\Lambda(h(x/N))}(1)}-\frac{1-\hat{\alpha}_v(\tfrac{\tilde{x}}{N})}{\hat{\alpha}_v(\tfrac{\tilde{x}}{N})}\right][\sqrt{f(\sigma^{x,v}\eta)}]^2 \, d\nu^N_{h}\\
				&-\displaystyle\frac{N^2}{4N^{\theta}}\int_{A_0}\sum_{\stackrel{x \in D^d_N}{x_1=1}}\sum_{v \in \mathcal{V}}\left(1-\hat{\alpha}_v\left(\tfrac{\tilde{x}}{N}\right)\right)\left[\frac{m_{\Lambda(h(x/N))}(0)}{1-m_{\Lambda(h(x/N))}(0)}-\frac{\hat{\alpha}_v(\tfrac{\tilde{x}}{N})}{1-\hat{\alpha}_v(\tfrac{\tilde{x}}{N})}\right][\sqrt{f(\sigma^{x,v}\eta)}]^2 \, d\nu^N_{h}.
			\end{split}
		\end{equation*}
		From  \eqref{eq:medidas} and since $f$ is a density, for  $\theta\in[0,1)$  the last display is equal to zero for $N$ sufficiently large, while for $\theta\geq 1$, it is bounded by 
		$
		\displaystyle\frac{C N^{d+1}}{N^{\theta}}.
		$
	\end{proof}
	
	\begin{remark}\label{obs:4}
		If $H(\mu^N|\nu^N_{h})$ is the relative entropy of the measure $\mu^N$ with respect to $\nu^N_{h}$, see \eqref{eq:1.12}, then there exists a constant $C_{h}$ such that $H(\mu^N|\nu^N_{h})\leq C_{h}N^d$. To prove it, note that by the definition of the entropy, 
		$$	H(\mu^N|\nu^N_{h}) =\displaystyle\int \log\left(\frac{\mu^N(\eta)}{\nu^N_{h}(\eta)}\right) \mu^N(\eta)\\[0.7cm]
		\leq \displaystyle\int \log\left(\frac{1}{\nu^N_{h}(\eta)}\right) \mu^N(\eta).	$$
		Since the measure $\nu^N_{h}$ is a product measure with marginal given by 
		$
		\nu^N_{h}\{\eta: \eta(x,\cdot)= \xi\}=m_{\Lambda(h(x/N))}(\xi),
		$
		where $m_{\lambda}(\cdot)$ was defined in \eqref{eq:1.2}, we obtain that the last display is bounded from above by
		\begin{equation}
			\begin{split}
				&\displaystyle\int \log\left(\frac{1}{\displaystyle\inf_{x \in D^d}(m_{\Lambda(h(x/N))}(1))\wedge(1-m_{\Lambda(h(x/N))}(1))}\right)^{N^d} \mu^N(\eta)\\
				=& N^d \log\left(\frac{1}{\displaystyle\inf_{x \in D^d}(m_{\Lambda(h(x/N))}(1))\wedge(1-m_{\Lambda(h(x/N))}(1))}\right).
			\end{split}
		\end{equation}			
		Since the functions $h_k^v$, defined in Remark \ref{obs:3}, are continuous, the image of each $h_k^v$ is a compact set bounded away from $0$ and $1$. Hence, from the definition of the measure $m$, we have $m_{\Lambda(h(x/N))}(1) > 0$ and $m_{\Lambda(h(x/N))}(1) < 1$. The constant can be taken as $$C_{h}:=\log\left(\displaystyle\frac{1}{\displaystyle\inf_{x \in D^d}m_{\Lambda(h(x/N))}(1)\wedge (1-m_{\Lambda(h(x/N))}(1))}\right).$$
	\end{remark}
	
	\subsection{Replacement Lemma at the Boundary}\label{subsec:6.2}
	Fix $k \in \{0,\dots, d\}$, a continuous function $G: [0,T]\times \mathbb{T}^{d-1} \to \mathbb{R}^{d+1}$, and for each $k$ consider the quantities
	\begin{equation*}
		\begin{split}
			&V_{k}^{1,\ell}(\eta_{s,\theta},  G) = \displaystyle\frac{1}{N^{d-1}}\displaystyle\sum_{\tilde{x} \in \mathbb{T}^{d-1}_N} G_k(s,\tfrac{\tilde{x}}{N})\left( I_k(\eta_{s,\theta}(1,\tilde{x}))-\sum_{v \in \mathcal{V}}v_k\alpha_v\left( \tfrac{\tilde{x}}{N}\right)\right),\\
			&V_{k}^{1,r}(\eta_{s,\theta},  G) = \displaystyle\frac{1}{N^{d-1}}\displaystyle\sum_{\tilde{x} \in \mathbb{T}^{d-1}_N} G_k(s,\tfrac{\tilde{x}}{N})\left( I_k(\eta_{s,\theta}(N-1,\tilde{x}))-\sum_{v \in \mathcal{V}}v_k\beta_v\left( \tfrac{\tilde{x}}{N}\right)\right),\\
			&V_{k}^{2,\ell}(\eta_{s,\theta}, G, \varepsilon) = \displaystyle\frac{1}{N^{d-1}}\displaystyle\sum_{\tilde{x} \in \mathbb{T}^{d-1}_N} G_k(s,\tfrac{\tilde{x}}{N})\Bigg( I_k(\eta_{s,\theta}(1,\tilde{x}))-\frac{1}{\lfloor \varepsilon N\rfloor}\sum_{x_1=1}^{\lfloor \varepsilon N\rfloor +1}I_k(\eta_{s,\theta}(x_1,\tilde{x}))\Bigg),\\
			&V_{k}^{2,r}(\eta_{s,\theta}, G, \varepsilon) = \displaystyle\frac{1}{N^{d-1}}\displaystyle\sum_{\tilde{x} \in \mathbb{T}^{d-1}_N} G_k(s,\tfrac{\tilde{x}}{N})\Bigg( I_k(\eta_{s,\theta}(N-1,\tilde{x}))-\frac{1}{\lfloor \varepsilon N\rfloor}\sum_{x_1=N-1-\lfloor \varepsilon N\rfloor}^{N-1}\!\!\!\!\!\!\!\!I_k(\eta_{s,\theta}(x_1,\tilde{x}))\Bigg),\\
		\end{split}
	\end{equation*}\noindent
	where $s \in [0,T]$, above $G_k$ are the components of function $G$ and $\lfloor r\rfloor$ the integer part of $r$.
	The main result of this section is the following proposition:
	\begin{proposition}[Replacement Lemma at the boundary]\label{prop:2}
		For each $0 \leq t \leq T$, $k \in \{0,\dots,d\}$ and $G \in C(D^d)$, it holds:
		\begin{itemize}\item [(1)] for $\theta\in[0,1)$:
			\begin{equation}
				\begin{array}{ll}
					\displaystyle\limsup_{N \to \infty}\mathbb{E}_{\mu^N}\left[\left|\int_0^t ds\, V^{1,\vartheta}_k(\eta_{s,\theta},G)\right|\right]=0
				\end{array}
			\end{equation}
			\item [(2)] for $\theta\geq 1$
			\begin{equation}
				\begin{array}{ll}
					\displaystyle\limsup_{\varepsilon \to 0} \limsup_{N \to \infty}\mathbb{E}_{\mu^N}\left[\left|\int_0^t ds\, V^{2,\vartheta}_k(\eta_{s,\theta},G,\varepsilon)\right|\right]=0.
				\end{array}
			\end{equation}
		\end{itemize}
		Above $\vartheta\in\{\ell, r\}$.
	\end{proposition} 
	
	\begin{proof}
		We present the proof for $V^{1,\vartheta}_k,$ but the other case is analogous.	By the entropy inequality and Jensen's inequality for any $A>0$ the expectation in the statement of the proposition is bounded from above by
		\begin{equation}\label{eq:1.24}
			\begin{array}{l}
				\displaystyle\frac{H(\mu^N|\nu^N_{h})}{AN^d}+\frac{1}{AN^d}\log \mathbb{E}_{\nu^N_{h}}\left[\exp\left\{\left|\int_0^t ds\, AN^d V^{1,\vartheta}_k(\eta_{s,\theta},G)\right|\right\}\right].
			\end{array}
		\end{equation}
		By Remark \ref{obs:4}, the leftmost term in last display  is bounded by $\tfrac{C_{h}}{A}$, so we only need to show that the rightmost term vanishes as $N \to \infty$. Since $e^{|x|} \leq e^x+e^{-x}$ and
		$$
		\displaystyle\limsup_{N \to \infty} N^{-d}\log\{a_N+b_N\} \leq \max\{\limsup_{N \to \infty} N^{-d}\log (a_N),\limsup_{N \to \infty} N^{-d}\log (b_N)\},
		$$
		we may remove the absolute value from the expression \eqref{eq:1.24}. By the Feynman-Kac formula, see for instance \cite{Kipnis/Landim}, \eqref{eq:1.24} is bounded from above by
		$$
		\begin{array}{l}
			\displaystyle\frac{C_{h}}{A}+ t \sup_{f}\left\{ \int V^{1,\vartheta}_k(\eta,  G)f(\eta) \, d\nu^N_{h} + \displaystyle\frac{\langle \mathcal{L}_{N,\theta}\sqrt{f}, \sqrt{f}\rangle_{\nu^N_{h}}}{AN^{d}}\right\},
		\end{array}
		$$\noindent
		where the supremum is taken over all densities $f$ with respect to $\nu^N_{h}$. 
		Since $\theta\in[0,1)$, from Corollaries 
		\ref{coro:1}, \ref{coro:2} and \ref{coro:3} the rightmost term in last display is bounded from above by
		\begin{equation}
			-\displaystyle\frac{1}{4 AN^{d-2}}D_{\nu^N_{h}}^\theta(\sqrt{f})
		\end{equation}
		As a consequence of the next lemma we can bound the integral term in the supremum above by
		$$
		\frac{C \varepsilon N^\theta}{ 4} +\displaystyle\frac{C'}{4\varepsilon N^{d-1}}D^{b,\theta}_{\nu^N_{h}}(\sqrt{f}),
		$$
		and $C,C'$ are constants. Making the choice $\varepsilon=\tfrac{AC'}{N}$,
		we bound the supremum above by
		$$
		\begin{array}{l}
			\displaystyle\frac{C_{h}}{A}+ t \Big(\frac{C C' A N^{\theta-1}}{4}+\tilde R_N\Big).
		\end{array}
		$$
		Finally taking $A=N^{\tfrac{1-\theta}{2}}$ and noting that $\tilde R_N$ vanishes for $N$ sufficiently large,
		the  proof ends. 
		
	\end{proof}
	\begin{lemma}\label{lema:3}Let $f$ be a density function with respect to $\nu^N_{h}$. For any $B>0$, for every $k \in \{0,\dots,d\}$, and every $G\in C([0,T]\times \mathbb{T}^{d-1})$, there exists a constant $\varepsilon>0$ such that
		$$
		\langle V^{1,\vartheta}_k(\eta, G) , f(\eta)\rangle_{\nu^N_{h}} \leq 
		\frac{C \varepsilon N^\theta}{ 4} +\displaystyle\frac{C'}{4\varepsilon N^{d-1}}D^{b,\theta}_{\nu^N_{h}}(\sqrt{f})+\tilde R_N,
		$$
		where 
		$\vartheta\in\{\ell, r\}$,  $C,C'$ are constants and  $\tilde R_N$ vanishes for $N$ sufficiently large.
	\end{lemma}
	
	\begin{proof}
		We prove for $\vartheta=\ell$, since for $\vartheta=r$ the proof is entirely analogous. First of all, note that since $G$ is continuous and its domain $[0,T]\times D^d$ is compact, it is enough to prove the result with $G=\bf{1}$. We note that, from \eqref{eq:sumalphabetahat},
		\begin{equation}
			\begin{split}
				\langle V^{1,\ell}_k(\eta, G) , f(\eta)\rangle_{\nu^N_{h}}=&\frac{1}{N^{d-1}}\sum_{\tilde{x} \in \mathbb{T}^{d-1}_N}\int f(\eta)\left( I_k(\eta_{s}(1,\tilde{x}))- \sum_{v \in \mathcal{V}}v_k \hat\alpha_v\left(\tfrac{\tilde{x}}{N}\right)\right)\, d\nu^N_{h}\\
				=&\frac{1}{N^{d-1}}\sum_{\tilde{x} \in \mathbb{T}^{d-1}_N}\int f(\eta)\left( I_k(\eta_{s}(1,\tilde{x}))- \sum_{v \in \mathcal{V}}v_k \hat{\alpha}_v\left(\tfrac{\tilde{x}}{N}\right)\right)\, d\nu^N_{h}
				\\=& \frac{1}{N^{d-1}}\sum_{\tilde{x} \in \mathbb{T}^{d-1}_N} \sum_{v \in \mathcal{V}}v_k \int f(\eta)\left(\eta(1,\tilde{x},v)- \hat{\alpha}_v\left(\tfrac{\tilde{x}}{N}\right)\right)\, d\nu^N_{h}.
			\end{split}
		\end{equation}
		By summing and subtracting an appropriate term, last term is equal to
		\begin{equation}
			\begin{split}\label{eq:1.25}
				&\displaystyle\frac{1}{2}\Big|\frac{1}{N^{d-1}}\sum_{\tilde{x} \in \mathbb{T}^{d-1}_N} \sum_{v \in \mathcal{V}}v_k \int \left(f(\eta)-f(\sigma^{x,v}\eta)\right)\left(\eta(1,\tilde{x},v)-\hat{\alpha}_v\left(\tfrac{\tilde{x}}{N}\right)\right)\, d\nu^N_{h}\Big|\\ 
				+&\displaystyle\frac{1}{2}\Big|\frac{1}{N^{d-1}}\sum_{\tilde{x} \in \mathbb{T}^{d-1}_N} \sum_{v \in \mathcal{V}}v_k \int \left(f(\eta)+f(\sigma^{x,v}\eta)\right)\left(\eta(1,\tilde{x},v)- \hat{\alpha}_v\left(\tfrac{\tilde{x}}{N}\right)\right)\, d\nu^N_{h}\Big|,
			\end{split}
		\end{equation}
		where $x=(1,\tilde x)$.
		Applying Young's inequality, $ab \leq \displaystyle\frac{\varepsilon N^\theta a^2}{2}+\frac{b^2}{2\varepsilon N^\theta}$, on the first term of last display, we can bound it from above by
		\begin{equation}
			\begin{split}\label{eq:1.26}
				&\displaystyle\frac{\varepsilon N^\theta}{4}\Big|\frac{1}{N^{d-1}}\sum_{\tilde{x} \in \mathbb{T}^{d-1}_N} \sum_{v \in \mathcal{V}}v_k \int \left(\sqrt{f(\eta)}+\sqrt{f(\sigma^{x,v}\eta)}\right)^2\frac{\left(\eta(1,\tilde{x},v)- \hat{\alpha}_v\left(\tfrac{\tilde{x}}{N}\right)\right)^2}{r^N_{x}(\eta,\alpha)}\, d\nu^N_{h}\Big|\\
				+&\displaystyle\frac{1}{4\varepsilon}\Big|\frac{1}{N^{d-1+\theta}}\sum_{\tilde{x} \in \mathbb{T}^{d-1}_N} \sum_{v \in \mathcal{V}}v_k \int \left(\sqrt{f(\eta)}-\sqrt{f(\sigma^{x,v}\eta)}\right)^2r^N_{x}(\eta,\alpha)\, d\nu^N_{h}\Big|
			\end{split}
		\end{equation}
		where $r^N_{x}(\eta,\alpha)$ was defined in \eqref{eq:taxas} and this holds for any $\varepsilon>0$.
		Since $\left(\eta(1,\tilde{x},v)- \hat{\alpha}_v\left(\tfrac{\tilde{x}}{N}\right)\right)^2 \leq 1$, $r^N_x$ is bounded away from zero and  one, then \eqref{eq:1.26} is bounded from above by 
		$$
		\begin{array}{l}
			\displaystyle\frac{C \varepsilon N^\theta}{ 4} +\displaystyle\frac{C'}{4\varepsilon N^{d-1}}D^{b,\theta}_{\nu^N_{h}}(\sqrt{f}),
		\end{array}
		$$
		with $C$ and $C'$ constants.
		Now, we analyze the second term of \eqref{eq:1.25}. Note that
		\begin{equation*}
			\begin{split}
				&\displaystyle\frac{1}{2}\Big|\frac{1}{N^{d-1}}\sum_{\tilde{x} \in \mathbb{T}^{d-1}_N} \sum_{v \in \mathcal{V}}v_k \int \left(f(\eta)+f(\sigma^{x,v}\eta)\right)\left(\eta(1,\tilde{x},v)- \hat{\alpha}_v\left(\tfrac{\tilde{x}}{N}\right)\right)\, d\nu^N_{h}\Big|\\ 
				=&\displaystyle\frac{1}{2}\Big|\frac{1}{N^{d-1}}\sum_{\tilde{x} \in \mathbb{T}^{d-1}_N} \sum_{v \in \mathcal{V}}v_k \int f(\eta)\left(\eta(1,\tilde{x},v)- \hat{\alpha}_v\left(\tfrac{\tilde{x}}{N}\right)\right)\, d\nu^N_{h}\Big|\\ 
				+&\displaystyle\frac{1}{2}\Big|\frac{1}{N^{d-1}}\sum_{\tilde{x} \in \mathbb{T}^{d-1}_N} \sum_{v \in \mathcal{V}}v_k \int f(\sigma^{x,v}\eta)\left(\eta(1,\tilde{x},v)- \hat{\alpha}_v\left(\tfrac{\tilde{x}}{N}\right)\right)\, d\nu^N_{h}\Big|.
			\end{split}
		\end{equation*}
		Recall the proof of Corollary \ref{coro:3}. By splitting the last two integrals over the sets $A_0=\{\eta:\eta((1,\tilde{x}),v)=0\}$ and $A_1=\{\eta:\eta((1,\tilde{x}),v)=1\}$, last display becomes equal to 
		\begin{equation*}
			\begin{split}
				&\displaystyle\frac{1}{2}\Big|\frac{1}{N^{d-1}}\sum_{\tilde{x} \in \mathbb{T}^{d-1}_N} \sum_{v \in \mathcal{V}}v_k \Big(\int_{A_0} f(0,\tilde\eta)\left(- \hat{\alpha}_v\left(\tfrac{\tilde{x}}{N}\right)\right)\, d\nu^N_{h}+\int_{A_1} f(1,\tilde\eta)\left(1- \hat{\alpha}_v\left(\tfrac{\tilde{x}}{N}\right)\right)\, d\nu^N_{h}\Big)\Big|\\ 
				+&\displaystyle\frac{1}{2}\Big|\frac{1}{N^{d-1}}\sum_{\tilde{x} \in \mathbb{T}^{d-1}_N} \sum_{v \in \mathcal{V}}v_k \Big(\int_{A_0} f(1,\tilde \eta)\left(- \hat{\alpha}_v\left(\tfrac{\tilde{x}}{N}\right)\right)\, d\nu^N_{h}+\int_{A_1} f(0,\tilde\eta)\left(1- \hat{\alpha}_v\left(\tfrac{\tilde{x}}{N}\right)\right)\, d\nu^N_{h}\Big)\Big|.
			\end{split}
		\end{equation*}
		Doing simple computations last display becomes equal to
		\begin{equation*}
			\tilde R_N:=\displaystyle\frac{1}{2N^{d-1}}\Big|\sum_{\tilde{x} \in \mathbb{T}^{d-1}_N} \sum_{v \in \mathcal{V}}v_k\sum_{\tilde \eta} \Big( f(0,\tilde\eta)+f(1,\tilde \eta)\Big) \left(m_{\Lambda (h(1/N,\tilde x/N))}(1)- \hat{\alpha}_v\left(\tfrac{\tilde{x}}{N}\right)\right)\nu^N_h(\tilde \eta)\Big|.
		\end{equation*}
		Using  the fact that $f$ is a density and that the function $h$ satisfies (1) in Remark \ref{obs:3} and from \eqref{eq:medidas}, together with the fact that the set $V $ is finite, we conclude that $\tilde R_N$ vanishes for $N$ sufficiently large. This ends the proof.
		
	\end{proof}
	
	\begin{lemma}\label{lema:4}
		For every $k \in \{0,\dots, d\}$, and every $G\in C([0,T]\times D^d)$, 
		\begin{equation*}
			\limsup_{\varepsilon \to 0}\limsup_{N \to \infty} \langle V^{2,\vartheta}_k( \eta, \zeta, G), f(\eta)\rangle_{\nu^N_{h}}=0,
		\end{equation*}
		where $\vartheta\in\{\ell, r\}$ and $h$ is constant.
	\end{lemma}
	
	\begin{proof}
		First of all, note that since $G$ is continuous and its domain $[0,T] \times D^d$ is compact, it is enough to prove the result with $G={\bf 1}$. We will only prove for $\vartheta=\ell$, since for $\vartheta=r$ the proof is entirely analogous. Observe that
		\begin{equation*}
			\begin{split}
				\langle V^{2,\ell}_k(\eta, \zeta, {\bf 1}), f(\eta)\rangle_{\nu^N_{h}} &	=\displaystyle\frac{1}{N^{d-1}} \sum_{\tilde{x} \in \mathbb{T}^{d-1}_N} \int f(\eta)\left( I_k(\eta(1,\tilde{x}))-\frac{1}{\varepsilon N}\sum_{x_1=1+1}^{\varepsilon N+1}I_k(\eta(x_1,\tilde{x}))\right) \, d\nu^N_{h}\\
				& =\displaystyle\frac{1}{N^{d-1}} \sum_{\tilde{x} \in \mathbb{T}^{d-1}_N}\sum_{v \in \mathcal{V}}v_k \int f(\eta)\left( \eta(1,\tilde{x},v)-\frac{1}{\varepsilon N}\sum_{x_1=1+1}^{\varepsilon N+1}\eta(x_1,\tilde{x},v)\right) \, d\nu^N_{h}\\
				& =\displaystyle\frac{1}{N^{d-1}} \sum_{\tilde{x} \in \mathbb{T}^{d-1}_N}\sum_{v \in \mathcal{V}}v_k \int f(\eta)\left(\frac{1}{\varepsilon N}\sum_{x_1=1+1}^{\varepsilon N+1}\left[ \eta(1,\tilde{x},v)-\eta(x_1,\tilde{x},v)\right]\right) \, d\nu^N_{h}.
			\end{split}
		\end{equation*}
		By writing the term $\displaystyle\frac{1}{\varepsilon N}\sum_{x_1=1+1}^{\varepsilon N+1}\left[\eta(1,\tilde{x},v)-\eta(x_1,\tilde{x},v)\right]$ as a telescopic sum, we obtain that the last term is equal to
		$$
		\frac{1}{N^{d-1}} \sum_{\tilde{x} \in \mathbb{T}^{d-1}_N}\sum_{v \in \mathcal{V}}v_k \int f(\eta)\left(\frac{1}{\varepsilon N}\sum_{x_1=1+1}^{\varepsilon N+1}\sum_{y=1}^{x_1-1}\left\{\eta(y,\tilde{x},v)-\eta(y+1,\tilde{x},v)\right\}\right) \, d\nu^N_{h}.
		$$
		Writing this sum as twice its half, performing change of variables and noting that $h$ is constant, we obtain that the last display is equal to
		\begin{equation}\label{eq:2parte}
			\displaystyle\frac{1}{N^{d-1}}\!\!\sum_{\tilde{x} \in \mathbb{T}^{d-1}_N}\sum_{v \in \mathcal{V}}v_k \frac{1}{2\varepsilon N}\sum_{x_1=1+1}^{\varepsilon N+1}\sum_{y=1}^{x_1-1} \int \left(f(\eta)-f(\eta^{y,y+1,v})\right) \left(\eta(y,\tilde{x},v)-\eta(y+1,\tilde{x},v)\right)\, d\nu^N_{h}.
		\end{equation}
		Rewriting $f(\eta)-f(\eta^{y,y+1,v})$ as $\left(\sqrt{f(\eta)}-\sqrt{f(\eta^{y,y+1,v})}\right)\left(\sqrt{f(\eta)}+\sqrt{f(\eta^{y,y+1,v})}\right)$ and using Young's inequality, for all $B>0$, we obtain that \eqref{eq:2parte} is bounded from above by
		\begin{equation*}
			\begin{split}
				&\displaystyle\frac{1}{N^{d-1}} \sum_{\tilde{x} \in \mathbb{T}^{d-1}_N}\sum_{v \in \mathcal{V}}v_k \frac{B}{2\varepsilon N}\sum_{x_1=1+1}^{\varepsilon N+1}\sum_{y=1}^{x_1-1} \int \left(\sqrt{f(\eta)}-\sqrt{f(\eta^{y,y+1,v})}\right)^2\, d\nu^N_{h}\\ 
				+&\displaystyle\frac{1}{N^{d-1}}\!\!\!\! \sum_{\tilde{x} \in \mathbb{T}^{d-1}_N}\!\!\sum_{v \in \mathcal{V}}v_k \frac{1}{2B\varepsilon N}\!\!\!\sum_{x_1=1+1}^{\varepsilon N+1}\!\!\sum_{y=1}^{x_1-1} \int\!\! \left(\sqrt{f(\eta)}+\sqrt{f(\eta^{y,y+1,v})}\right)^2\!\! \left(\eta(y,\tilde{x},v)-\eta(y+1,\tilde{x},v)\right)^2 d\nu^N_{h}.
			\end{split}
		\end{equation*}
		Using that $f$ is a density for $\nu^N_{h}$, the second term in last display is bounded by $\frac{C\varepsilon N}{B}$. Letting the sum in $y$ run from $1$ to $N-1$, the first term in last display is bounded by $BD^{ex}_{\nu^N_{h}}(\sqrt{f})$. By Corollary \ref{coro:1}, since $h$ is a constant function, we obtain
		$$
		D^{ex}_{\nu^N_{h}}(\sqrt{f})=-\langle \mathcal{L}_N^{ex} \sqrt{f}, \sqrt{f}\rangle_{\nu^N_{h}}.
		$$
		Since $0 \leq D^c_{\nu^N_{h}}(\sqrt{f}) = -\langle \mathcal{L}_N^{c} \sqrt{f}, \sqrt{f}\rangle_{\nu^N_{h}}$ and $0 \leq D^b_{\nu^N_{h}}(\sqrt{f})$, using Corollary \ref{coro:3}, we have that $D^{b,\theta}_{\nu^N_{h}}(\sqrt{f})= -\langle \mathcal{L}_{N,\theta}^{b} \sqrt{f}, \sqrt{f}\rangle_{\nu^N_{h}}$.
		Therefore,
		\begin{equation}\nonumber
			D^{ex}_{\nu^N_{h}}(\sqrt{f}) \leq - B\langle \mathcal{L}_{N,\theta} \sqrt{f}, \sqrt{f}\rangle_{\nu^N_{h}}.
		\end{equation}
	\end{proof}
	
	\subsection{Replacement lemma at the bulk}
	Before we state the replacement lemma at the bulk that will allow us to prove that the limit points $\mathbb{Q}^*_\theta$ are concentrated on weak solutions of the system of partial differential equations \eqref{eq:PDE}, we introduce some notations. Recall the definition of $\textbf{I}^L$ in \eqref{eq:def_boldI}. 
	
	\begin{lemma}[Replacement lemma at the bulk]\label{lema:5}
		For all $\delta >0$,  $j \in\{1,\dots,d\}$ and $k \in \{0,\dots,d\}$:
		\begin{equation*}
			\limsup_{\varepsilon \to 0}\limsup_{N \to \infty}\mathbb{P}_{\mu^N}\left[  \int_0^T \frac{1}{N^d} \sum_{x \in D^d_N} \tau_xV_{\varepsilon N}^{j,k}(\eta_{s,\theta})\, ds \geq \delta\right]=0,
		\end{equation*}
		where 
		\begin{equation*}
			V_{\ell}^{j,k}(\eta)=\Big|\displaystyle \frac{1}{(2\ell+1)^d}\sum_{y \in \Lambda_{\ell}}\sum_{v \in \mathcal{V}}v_k\sum_{z \in \mathbb{Z}^d}p(z,v)z_j\tau_y[\eta(0,v)(1-\eta(z,v))]-\displaystyle\sum_{v \in \mathcal{V}}v_jv_k\chi(\theta_v(\Lambda(\textbf{I}^{\ell}(0,\eta))))\Big|.
		\end{equation*}
	\end{lemma}
	Note that for all $j \in \{1,\dots,d\}$ and $k \in \{0,1,\dots,d\}$, $V^{j,k}_{\varepsilon N}$ is well-defined for large $N$ since $p(\cdot, v)$ is of finite range. We now observe that Corollaries \ref{coro:1} and \ref{coro:2} permit us to prove the previous replacement lemma for the boundary driven exclusion process by using the process without the boundary part of the generator (see \cite{Kipnis/Landim} for further details). For the proof of Lemma \ref{lema:5}, see \cite{Simas}. 
	
	\section{Energy Estimates}\label{sec:7}
	We will now define some quantities in order to prove that each component of the vector solution belongs, in fact, to $L^2([0,T],\mathscr{H})$.
	Let the energy $\mathscr{E}: \mathcal{D}([0,T], \mathcal{M}) \to [0,\infty]$ be given by
	\begin{equation*}
		\mathscr{E}(\pi)=\displaystyle\sum_{i=1}^d \mathscr{E}_i(\pi),
	\end{equation*}
	with 
	$	\mathscr{E}_i(\pi)=\displaystyle\sup_{G \in C^{\infty}_{c}(\Omega_T)}\left\{ 2\int_0^T\,dt\, \langle \pi_t,\partial_{u_i}G_t\rangle - \int_0^T \, dt \int_{D^d} \,du \,G(t,u)^2\right\},
	$
	where $\Omega_T=(0,T)\times D^d$ and $C^\infty_c(\Omega_T)$ stands for the set of infinitely differentiable functions (with respect to time and space) with compact support contained in $\Omega_T$. For any $G \in C^{\infty}_c(\Omega_T), \, 1 \leq i \leq d \mbox{ and } C\geq 0$, let the functional $\mathscr{E}_{i,C}^G: \mathcal{D}([0,T], \mathcal{M}) \to \mathbb{R}$ be given by
	\begin{equation*}
		\mathscr{E}_{i,C}^G(\pi)=\int_0^T ds\,   \langle \pi_s,\partial_{u_i}G_s\rangle - C\int_0^T  ds \int_{D^d} du \,G(s,u)^2.
	\end{equation*} 
	Note that
	\begin{equation}\label{eq:1.36}
		\sup_{G \in C^\infty_c(\Omega_T)}\{\mathscr{E}_{i,C}^G\}= \displaystyle\frac{\mathscr{E}_i(\pi)}{4C}.
	\end{equation}
	It is well-known that $\mathscr{E}(\pi)$ is finite if, and only if, $\pi$ has a generalized gradient, $\nabla \pi = (\partial_{u_1}\pi, \dots , \partial_{u_d}\pi)$, which is a measurable function and 
	\begin{equation*}
		\tilde{\mathscr{E}}(\pi)= \int_0^T ds\, \int_{D^d} du \, \|\nabla\pi_t(u)\|^2 < \infty,
	\end{equation*}
	in which case, $\mathscr{E}(\pi)=\tilde{\mathscr{E}}(\pi)$. Recall from Section \ref{sec:4} that the sequence $(\mathbb{Q}_{N,\theta})_N$ is tight. We have that: 
	
	\begin{proposition}\label{prop:1}
		Let $\mathbb{Q}^*_\theta$ be any limit point of the sequence of measures $(\mathbb{Q}_{N,\theta})_N$. Then, 
		\begin{equation}\nonumber
			E_{\mathbb{Q}^*_\theta}\left[\int_0^Tds\, \left(\int_{D^d}\| \nabla \varrho^k(s,u)\|^2 \, du\right)\right] < \infty,
		\end{equation}
		for $k \in \{0, 1,\dots,d\}.$
	\end{proposition}
	
	The proof of Proposition \ref{prop:1} is given after the following lemmas.
	\begin{lemma}\label{lema:6}
		For all $\theta \geq 0$ and for each $i \in \{1,\dots,d\}$ there is a positive constant $C>0$ such that 
		\begin{equation*}
			E_{\mathbb{Q}_\theta^{*}}\left[\sup_{G}\left\{ \int_0^T\int_{D^d} \partial_{u_i}G(s,u)\varrho^k(s,u)du\, ds - C\int_0^T ds \int_{D^d} du\, G(s,u)^2\right\}\right] < \infty,
		\end{equation*}
		for $k \in \{0,1,\dots,d\}$, where the supremum is carried over all the functions $G \in C^{\infty}(\Omega_T)$ and $\varrho^0=\rho.$
	\end{lemma}
	
	\begin{proof}
		Let $\{G^m :  m \geq 1 \}$ be a sequence of functions in $C^{\infty}_c(\Omega_T)$ dense in $L^2(\Omega_T)$. Thus, it is sufficient to prove that, for every $r\geq 1$,
		\begin{equation}\label{eq:1.37}
			E_{\mathbb{Q}^\ast_\theta}\left[\max_{1 \leq m \leq r}\left\{ \mathscr{E}^{G^m}_{i,C}(\pi^{k})\right\}\right]\leq \tilde{C},
		\end{equation}
		for some constant $\tilde{C}>0$, independent of $r$. 
		The expression on the left-hand side of \eqref{eq:1.37} is equal to
		\begin{equation}\label{eq:1.38}
			\lim_{N \to \infty}E_{\mu^N}\left[\max_{1 \leq m \leq r}\left\{ \int_0^T\langle \partial_{u_i}G^m(s,u),\pi^{k,N}_{s,\theta}\rangle ds - C\int_0^T ds \int_{D^d} du\, G^m(s,u)^2\right\}\right].
		\end{equation}
		By the relative entropy bound (see Remark \ref{obs:4}), Jensen's inequality and the fact that $\exp\{\displaystyle\max_{1 \leq j \leq k} a_j \} \leq \sum_{1 \leq j\leq k}\exp{a_j}$, the expectation in \eqref{eq:1.38} is bounded from above by
		\begin{equation*}
			\displaystyle\frac{H(\mu^N|\nu^N_{h})}{N^d}
			+ \displaystyle\frac{1}{N^d}\log\sum_{1 \leq m \leq r}E_{\nu^N_{h}}\left[\exp\left\{\int_0^T N \langle \partial_{u_i}G^m(s,u),\pi^{k,N}_{s,\theta}\rangle ds - C\int_0^T ds \int_{D^d} du\, G^m(s,u)^2\right\}\right],
		\end{equation*}
		where the functions $h$ are the same as in Remark \ref{obs:3} in Section \ref{sec:5}.
		
		We can bound the first term in the sum above by $C_h$. It is enough to show, for a fixed function $G$, that
		\begin{equation*}
			\displaystyle\limsup_{N \to \infty}\frac{1}{N^d}\log E_{\nu^N_{h}}\left[\exp\left\{\int_0^T N \langle \partial_{u_i}G(s,u),\pi^{k,N}_{s,\theta}\rangle ds - C\int_0^T ds \int_{D^d} du\, G(s,u)^2\right\}\right] \leq \tilde{c}
		\end{equation*}
		for some constant $\tilde{c}$ independent of $G$. 
		Then the result follows from the next lemma and the definition of the empirical measure.
	\end{proof}
	
	\begin{lemma}\label{lema:7}
		There exists a constant $C_0=C_0(h)>0$, such that for every $i=1,\dots,d$ every $k \in \{0,\dots,d\}$ and every function $G \in C^{\infty}_c(\Omega_T)$
		\begin{equation*}
			\limsup_{N \to \infty}\frac{1}{N^d}\log E_{\nu^N_{h}}\left[\exp \{ N^d \mathscr{E}_{i,C_0}^{G}(\pi^{k,N}_\theta)\}\right] \leq C_0.
		\end{equation*}
	\end{lemma}
	
	\begin{proof}
		Writing $\partial_{u_i}G_s\left(\tfrac{x}{N}\right)=N\left[G_s\left(\tfrac{x+e_i}{N}\right)-G_s\left(\tfrac{x}{N}\right)\right]+O(N^{-1})$ and summing by parts (the compact support of $G$ takes care of the boundary term), by applying the Feynman-Kac formula and using the same arguments as in the proof of Lemma \ref{prop:2}, we have that 
		\begin{equation*}
			\displaystyle\frac{1}{N^d}\log E_{\nu^N_{h}}\left[\exp \big\{ N \int_0^T \, ds \sum_{x \in D^d_N}\left(I_k(\eta_{s,\theta}(x))-I_k(\eta_{s,\theta}(x-e_i))\right)G\left(s,\tfrac{x}{N}\right)\big\}\right] \leq  
			\displaystyle\frac{1}{N^d}\int_0^T\lambda_s^N\, ds,
		\end{equation*}
		where $\lambda^N_s$ is equal to 
		\begin{equation}\label{eq:1.39}
			\sup_{f}\left\{ \left\langle N \sum_{x \in D^d_N}((I_k(\eta_{s,\theta}(x))-I_k(\eta_{s,\theta}(x-e_i)))G\left(s,\tfrac{x}{N}\right), f\right\rangle_{\nu^N_{h}}+N^2\langle \mathcal{L}_{N,\theta}\sqrt{f}, \sqrt{f}\rangle_{\nu^N_{h}}\right\},
		\end{equation}
		where the supremum is taken over all densities $f$ with respect to $\nu^N_{h}$. 
		Now we will consider two cases: if $h$ is a constant function, then from Corollaries \ref{coro:1}, \ref{coro:2} and \ref{coro:3}, the expression inside brackets is bounded from above by
		\begin{equation*}
			-\frac{N^2}{2}D_{\nu^N_{h}}(\sqrt{f})+\sum_{x \in D^d_N}\left\{ NG\left(s,\tfrac{x}{N}\right)\int[I_k(\eta_{s,\theta}(x))-I_k(\eta_{s,\theta}(x-e_i))]f(\eta)d\nu^N_{h}\right\}.
		\end{equation*}
		We now rewrite the term inside the brackets as 
		\begin{equation}\label{eq:1.40}
			\sum_{v \in \mathcal{V}}v_k\sum_{x \in D^d_N}\left\{ \int NG\left(s,\tfrac{x}{N}\right)[\eta(x,v)-\eta(x-e_i,v)]f(\eta)d\nu^N_{h}\right\}.
		\end{equation}
		After a simple computation, we may rewrite the terms inside the brackets of the above expression as
		\begin{equation}\label{eq:limitad}
			\begin{split}
				&	NG\left(s,\tfrac{x}{N}\right)\displaystyle\int[\eta(x,v)-\eta(x-e_i,v)]f(\eta)d\nu^N_{h}\\
				= & NG\left(s,\tfrac{x}{N}\right)\displaystyle\int \eta(x,v)f(\eta)\,d\nu^N_{h}
				-NG\left(s,\tfrac{x}{N}\right)\displaystyle\int \eta(x,v)f(\eta^{x-e_i,x,v})\frac{\nu^N_{h}(\eta^{x,x-e_i,v})}{\nu^N_{h}(\eta)}d\nu^N_{h} \\
				= & NG\left(s,\tfrac{x}{N}\right)\displaystyle\int \eta(x,v)[f(\eta)-f(\eta^{x-e_i,x,v})]d\nu^N_{h}.
			\end{split}
		\end{equation}
		By using  $f(\eta)-f(\eta^{x-e_i,x,v})=\left(\sqrt{f(\eta)}-\sqrt{f(\eta^{x-e_i,x,v})}\right)\left(\sqrt{f(\eta)}+\sqrt{f(\eta^{x-e_i,x,v})}\right)$ and applying Young's inequality, the equation \eqref{eq:limitad} is bounded from above by 
		$$
		\displaystyle\frac{N^2}{2} \displaystyle\int [\sqrt{f(\eta^{x-e_i,x,v})}-\sqrt{f(\eta)}]^2d\nu^N_{h} + 2G\left(s,\tfrac{x}{N}\right)^2\displaystyle\int\eta(x,v)(\sqrt{f(\eta)}+\sqrt{f(\eta^{x-e_i,x,v})})^2d\nu^N_{h}.\\
		$$
		Using the above estimate \eqref{eq:1.40} is clearly bounded by $\tfrac{N^2}{2} D_{\nu^N_{h}}(\sqrt{f})+CG\left(s,\tfrac{x}{N}\right)^2$, where $C$ is a positive constant. Thus, letting $C_0=C$, the statement of the lemma holds.
		
		Now we will analyze the case when $h$ is a smooth function. By Corollaries \ref{coro:1}, \ref{coro:2} and \ref{coro:3}, the expression inside brackets in \eqref{eq:1.39} is bounded from above by
		\begin{equation}\nonumber
			CN^d -\frac{N^2}{4}D_{\nu^N_{h}}(\sqrt{f})+\sum_{x \in D^d_N}\left\{ NG\left(s,\tfrac{x}{N}\right)\int[I_k(\eta_x(s))-I_k(\eta_{x-e_i}(s))]f(\eta)d\nu^N_{h}\right\}.
		\end{equation}
		Rewriting the term above, we will analyze the expression 
		\begin{equation}\label{eq:1.41}
			\sum_{x \in D^d_N}\sum_{v \in \mathcal{V}}v_k\left\{ NG\left(s,\tfrac{x}{N}\right)\int[\eta(x,v)-\eta(x-e_i,v)]f(\eta)d\nu^N_{h}\right\}.
		\end{equation}
		Now rewrite the term inside the brackets as 
		\begin{equation*}
			\begin{split}
				&NG\left(s,\tfrac{x}{N}\right)\displaystyle\int[\eta(x,v)-\eta(x-e_i,v)]f(\eta)d\nu^N_{h}\\ 
				=&NG\left(s,\tfrac{x}{N}\right)\displaystyle\int \eta(x,v)f(\eta)d\nu^N_{h}
				-NG\left(s,\tfrac{x}{N}\right)\displaystyle\int \eta(x,v)f(\eta^{x-e_i,x,v})\tfrac{\nu^N_{h}(\eta^{x,x-e_i,v})}{\nu^N_{h}(\eta)}d\nu^N_{h} \\ 
				=&NG\left(s,\tfrac{x}{N}\right)\displaystyle\int \eta(x,v)[f(\eta)-f(\eta^{x-e_i,x,v})]d\nu^N_{h}\\
				+&G\left(s,\tfrac{x}{N}\right)\displaystyle\int \eta(x,v)f(\eta^{x-e_i,x,v})N\left[1-\tfrac{\nu^N_{h}(\eta^{x,x-e_i,v})}{\nu^N_{h}(\eta)}\right]d\nu^N_{h}.
			\end{split}
		\end{equation*}
		Since $f(\eta)-f(\eta^{x-e_i,x,v})=\left(\sqrt{f(\eta)}-\sqrt{f(\eta^{x-e_i,x,v})}\right)\left(\sqrt{f(\eta)}+\sqrt{f(\eta^{x-e_i,x,v})}\right)$ and applying Young's inequality, the expression is bounded from above by 
		\begin{equation*}
			\begin{split}
				&N^2\displaystyle\int \frac{1}{2}[\sqrt{f(\eta^{x-e_i,x,v})}-\sqrt{f(\eta)}]^2d\nu^N_{h} 
				+ 2G\left(s,\tfrac{x}{N}\right)^2\displaystyle\int\eta(x,v)(\sqrt{f(\eta)}+\sqrt{f(\eta^{x-e_i,x,v})})^2d\nu^N_{h}\\ 
				+ & G\left(s,\tfrac{x}{N}\right)^2\displaystyle\int f(\eta^{x-e_i,x,v})d\nu^N_{h}
				+ \displaystyle\frac{1}{4}\int \eta(x,v)f(\eta^{x-e_i,x,v})\left[ N\left( 1-\frac{\nu^N_{h}(\eta^{x,x-e_i,v})}{\nu^N_{h}(\eta)}\right)\right]^2 d\nu^N_{h}.
			\end{split}
		\end{equation*}
		Using the above estimate, \eqref{eq:1.41} is clearly bounded by $C_1+C_1G\left(s,\tfrac{x}{N}\right)^2$, by some positive constant $C_1=C_1(h)$, using the estimate \eqref{eq:1.18} and the fact that $f$ is a density with respect to $\nu^N_{h}$. Thus, letting $C_0=C+C_1$, the statement of the lemma follows. 
	\end{proof}
	
	\begin{proof}[Proof of Proposition \ref{prop:1}]\label{proof}
		Let $\{G_m : 1 \leq m \leq r\}$ be a sequence of functions in $C^{\infty}_c(\Omega_T)$ dense in $L^2(\Omega_T)$. By the entropy inequality, see Remark \ref{obs:4}, there exists a constant $C_{h}>0$ such that 
		\begin{equation*}
			E_{\mu^N}\left[ \max_{1 \leq m \leq r }\Big\{\mathscr{E}^{G_m}_{i,C_0}(\pi^{k,N})\Big\}\right] \leq C_{h} + \frac{1}{N^d}\log E_{\nu^N_{h}}\left[\exp\Big\{ N^d\max_{1 \leq m \leq r} \{ \mathscr{E}^{G_m}_{i,C_0}(\pi^{k,N})\}\Big\}\right],
		\end{equation*}
		with $ i \in \{1,\dots,d\}$, and $k \in \{0,\dots,d\}$. Therefore, using Lemma \ref{lema:7} together with the elementary inequalities
		\begin{equation*}
			\displaystyle\limsup_{N \to \infty} N^{-d}\log(a_N+b_N) \leq \limsup_{N \to \infty}\max\Big\{\limsup_{N \to \infty}N^{-d}\log(a_N),\limsup_{N \to \infty}N^{-d}\log(b_N)\Big\}
		\end{equation*}
		and
		\begin{equation*}
			\exp\{\max\{x_1,\dots,x_n\}\} \leq \exp(x_1)+\dots +\exp(x_n)
		\end{equation*}
		we set that
		\begin{equation*}
			E_{\mathbb{Q}^*}\left[ \displaystyle\max_{1 \leq m \leq r }\Big\{\mathscr{E}^{G_m}_{i,C_0}(\pi^{k,N})\Big\}\right] = \displaystyle\lim_{N \to \infty} E_{\mu^N}\left[ \max_{1 \leq m \leq r }\Big\{\mathscr{E}^{G_m}_{i,C_0}(\pi^{k,N})\Big\}\right] \leq C_{h}+C_0.
		\end{equation*}
		Using this, the equation \eqref{eq:1.36} and the monotone convergence Theorem, we obtain the desired result.
	\end{proof}
	
	\appendix
	\section{Uniqueness of Weak Solutions}\label{appendix}
	To conclude the proof of the hydrodynamic limit, it remains to prove the uniqueness of weak solutions to \eqref{eq:1.6}. In the following, consider $\mathfrak{C}=+\infty$ and $(\rho^1,\varrho^1)$, $(\rho^2,\varrho^2)$ two weak solutions of \eqref{eq:1.6} with the same initial condition. Denote their difference by $(\bar{\rho},\bar{\varrho})=(\rho^1-\rho^2,\varrho^1-\varrho^2)$.  Note that $(\bar{\rho},\bar{\varrho})=\bar{\varrho}^k=0$ if, and only if, each component is equal to zero, which means that $\bar{\varrho}^k=0$ for each $k\in \{0,1,\dots,d\}$. Let $\{\psi_z\}_z$ be given by $\psi_z(u)=\sqrt{2}\sin(\pi z\,\cdot \,u)$ for $\|z\|\geq 1$ and $\psi_0(u)=1$ which is an orthonormal basis of $L^2(D^d)$. Let
	\begin{equation*}
		V_k(t)=\displaystyle\sum_{z\geq 0}\frac{1}{2a_z}\langle \bar{\varrho}^k_t , \psi_z \rangle^2
	\end{equation*}
	where $a_z=\|z\|^2\pi^2+1$. We claim that $V_k'(t) \leq CV_k(t)$, where $C$ is a positive constant. Since $V_k(0)=0, \, \forall \, k \in \{0, 1,\dots,d\}$, from Gronwall's inequality we will conclude that $V_k(t)\leq 0$, but since we know by definition that $V_k(t)\geq 0 $, we are done. Now we need to show that the claim is true.
	Note that 
	\begin{equation*}
		V'_k(t)=\displaystyle\sum_{z\geq 0}\frac{1}{a_z}\langle \bar{\varrho}^k_t, \psi_z \rangle \frac{d}{dt}\langle \bar{\varrho}^k_t, \psi_z \rangle,
	\end{equation*}
	and from the integral formulation \eqref{eq:1.6} we have that 
	\begin{equation*}
		\begin{split}
			\frac{d}{dt}\langle \bar{\varrho}^k_t, \psi_z \rangle = &  \langle \frac{d}{dt}\bar{\varrho}^k_t , \psi_z \rangle  +\langle \bar{\varrho}^k_t, \frac{d}{dt}\psi_z \rangle\\
			= &\frac{1}{2} \langle \bar{\varrho}^k_t, \Delta \psi_z \rangle +\langle \chi(\theta_v(\Lambda(\rho^1_t,\varrho^1_t)))-\chi(\theta_v(\Lambda(\rho^2_t,\varrho^2_t))), \partial_u \psi_z \rangle. 
		\end{split}
	\end{equation*}
	Since $\psi_z(u)=\sqrt{2}\sin(\pi z\,\cdot \, u)$, we have that $\partial_u\psi_z(u)=\sqrt{2}\pi \|z\|\cos(\pi z\, \cdot \, u)$ and $\Delta\psi_z(u)=-(\|z\|\pi)^2\sqrt{2}\sin(\pi z \, \cdot \, u)=-(z\pi)^2\psi_z$, then
	\begin{equation*}
		V'_k(t)=\displaystyle\sum_{z\geq 0}\frac{-(\|z\|\pi)^2}{2a_z} \langle \bar{\varrho}^k_t, \psi_z \rangle^2 + \sum_{z\geq 0}\frac{1}{a_z}  \langle \bar{\varrho}^k_t, \psi_z \rangle\langle \chi(\theta_v(\Lambda(\rho^1_t,\varrho^1_t)))-\chi(\theta_v(\Lambda(\rho^2_t,\varrho^2_t))), \partial_u \psi_z \rangle .
	\end{equation*}
	Using Young's inequality on the second term on the right-hand side of last identity, we bound that term from above by
	$$
	\displaystyle\frac{1}{2A}\sum_{z \geq 0}\frac{1}{a_z}\langle \bar{\varrho}^k_t, \psi_z \rangle^2+\frac{A}{2}\sum_{z \geq 0}\frac{1}{a_z}\langle \chi(\theta_v(\Lambda(\rho_t^1,\varrho_t^1)))-\chi(\theta_v(\Lambda(\rho_t^2,\varrho_t^2))),\partial_u\psi_z \rangle^2, \forall A>0.
	$$
	Observe that $\partial_u\psi_z=\|z\|\pi \phi_z(u)$, with $\phi_z(u)=\sqrt{2}\cos(\pi z \, \cdot\, u)$ for $\|z\| \geq 1$ and $\phi_0(u)=1$. Therefore, the second term at right-hand side in last display can be rewritten as
	\begin{equation*}
		\begin{split}
			&\displaystyle\frac{A}{2}\sum_{z \geq 0}\frac{(\|z\|\pi)^2}{a_z}\langle \chi(\theta_v(\Lambda(\rho_t^1,\varrho_t^1)))-\chi(\theta_v(\Lambda(\rho_t^2,\varrho_t^2))),\phi_z \rangle^2\\
			\leq &\displaystyle\frac{A}{2}\sum_{z \geq 0}\langle \chi(\theta_v(\Lambda(\rho_t^1,\varrho_t^1)))-\chi(\theta_v(\Lambda(\rho_t^2,\varrho_t^2))),\phi_z \rangle^2
		\end{split}
	\end{equation*}
	because of the choice for $a_z$. Observe that, since $\{\phi_z\}_z$ is an orthonormal basis of $L^2(D^d)$ we can rewritten the last display as 
	\begin{equation*}
		\displaystyle\frac{A}{2}\int_0^1 \left( \chi(\theta_v(\Lambda(\rho_t^1,\varrho_t^1)))-\chi(\theta_v(\Lambda(\rho_t^2,\varrho_t^2))) \right)^2 du.
	\end{equation*}
	Since $\chi(\theta_v(\Lambda(\cdot)))$ is Lipschitz, the last display is bounded from above by $\tfrac{A}{2}\|\bar{\varrho}_t\|^2_2$. Putting all this together, we conclude that
	\begin{equation*}
		V'_k(t)\leq \displaystyle\sum_{z \geq 0}\left( \frac{-(\|z\|\pi)^2}{2a_z}+\frac{1}{2Aa_z}+\frac{A}{2}\right)\langle \bar{\varrho}^k_t, \psi_z \rangle^2.
	\end{equation*}
	Taking $A=1$, we get 
	\begin{equation*}
		V'_k(t)  \leq \displaystyle\sum_{z \geq 0}\left( \frac{1}{2a_z}+\frac{1}{2}\right)\langle \bar{\varrho}^k_t, \psi_z \rangle^2  = \displaystyle\frac{1+a_z}{2a_z}\langle \bar{\varrho}^k_t, \psi_z \rangle^2 = C \,V_k(t).
	\end{equation*}
	The proof of uniqueness of weak solutions of \eqref{eq:1.6} with $\mathfrak{C}=0$ is similar to the previous one, considering the orthonormal basis $\phi_z(u)=\sqrt{2}\cos(\pi z\, \cdot \, u)$ of $L^2(D^d)$. Details are omitted here. We note that adapting the proof above to the case $\mathfrak{C}=1$ does no longer work since we could use the orthonormal basis of $L^2(D^d)$ given by a linear combination of $\sin$'s and $\cos$'s, but the problem is that, when we derive this basis the resulting set is no longer orthonormal and so the argument does not follow. Fortunately, we have an answer about the uniqueness of the $1$-dimensional case, it is possible to adapt the strategy used in \cite{Franco/Goncalves/Landim/Neumann} for our case. For higher dimensions the proof is left open.
	
	\section{Auxiliary results}\label{app:diffeo}
	In this Appendix we show that under Hypothesis \ref{eq:conjV}, the map $(\rho,\varrho)$ is a diffeomorphism onto its image.
	
	Begin by observing that the map $(\rho,\varrho)$ is the gradient of the function $F(\lambda) := \log(Z(\lambda))$, where $Z(\lambda)$ is the normalizing constant given in \eqref{eq:1.2}. The result, thus, follows by showing that $F(\cdot)$ has a strictly positive-definite Hessian. Indeed, since the Hessian of $F(\cdot)$ is the derivative of $\nabla F(\cdot)$, this also shows that $\nabla F(\cdot)$ has a non-vanishing derivative, and, thus, by the inverse function theorem, $(\rho,\varrho)(\cdot) = \nabla F(\cdot)$ is a local diffeomorphism. Further, this also implies that $F(\cdot)$ is a strictly convex function, and it is well-known that the gradient of a strictly convex function on a convex domain is injective, which shows that $(\rho,\varrho)(\cdot) = \nabla F(\cdot)$ is, actually, a global diffeomorphism onto its image.

	We now show that the Hessian of $F$, namely $H(\lambda)$, is strictly positive-definite. To this end, let $\alpha\in\mathbb{R}^{d+1}$ and observe that
	\begin{equation*}
		\begin{split}
			\langle \alpha, H(\lambda) \alpha\rangle &= \sum_{k,j=1}^{d+1} \sum_{v\in\mathcal{V}} \alpha_k \vartheta_k(v) \vartheta_j(v) \alpha_j\\
			&= \sum_{v\in \mathcal{V}} \left(\langle \alpha, \vartheta(v)\rangle \right)^2,
		\end{split}
	\end{equation*}
	where $\langle\cdot,\cdot\rangle$ denotes the standard inner product in $\mathbb{R}^{d+1}$, $\vartheta(v) := (\vartheta_1(v),\ldots, \vartheta_{d+1}(v))$, and $\vartheta_j(v) = v_{j-1} \exp\{0.5\langle \lambda, \widetilde{v}\rangle\}/(1 + \exp\{\langle \lambda, \widetilde{v}\rangle\}), j=1,\ldots, d+1,  \widetilde{v} = (1,v), v\in\mathcal{V}$ and we define $v_0=1$. Therefore, 
	$$\langle \alpha, H(\lambda) \alpha\rangle = 0 \Longleftrightarrow \forall v\in\mathcal{V}, \alpha \perp \vartheta(v).$$
	Furthermore, this means that 
	$$\forall v\in\mathcal{V}, \sum_{k=1}^{d+1} \alpha_k v_{k-1} \exp\{0.5\langle \lambda, \widetilde{v}\rangle\}/(1 + \exp\{\langle \lambda, \widetilde{v}\rangle\}) = 0,$$
	so that, for any $v\in\mathcal{V}$, 
	$$\alpha \perp \vartheta(v) \Longleftrightarrow \alpha \perp \widetilde{v}.$$
	Thus, the Hessian, $H(\cdot)$ will be strictly positive-definite if, and only if, the following condition holds
	$$\alpha \perp \widetilde{v} \text{ for all $v\in\mathcal{V}$} \Rightarrow \alpha = 0.$$
	Observe that all we need to hold is the above condition. We now show that Hypothesis \ref{hyp_1} is sufficient. To this end, let $\alpha\in\mathbb{R}^{d+1}$ be such that 
	$$\forall v\in\mathcal{V}, \langle \alpha, \widetilde{v}\rangle = 0.$$
	Take any $v\in\mathcal{V}$ such that there exists $j$ such that $v_j\neq 0$. From the invariance by permutations,  there exists some $v\in\mathcal{V}$ such that $v_1\neq 0$. Thus,
	$$\alpha_{2} = -\alpha_1/v_1 - \alpha_3 v_2/v_1 - \cdots - \alpha_{d+1} v_{d}/v_1.$$
	From the  invariance by reflections, we obtain that 
	$$\alpha_2 = -\alpha_2 \Rightarrow \alpha_2 = 0.$$
	By repeating the same argument for $v_2,\ldots,v_d$, using the invariance by permutations to ensure that for each $j=2,\ldots,d$, there exists $v\in\mathcal{V}$, such that $v_j\neq 0$, we obtain that
	$$\alpha_2 = \alpha_3 = \cdots = \alpha_{d+1} = 0.$$
	Finally, we use $\langle \alpha, \widetilde{v}\rangle = 0$ one last time to conclude that $\alpha_1 = 0$, which shows that $\alpha = 0$ and concludes the proof.
	\\
	
	\section*{Acknowledgments}\noindent
	O.A. thanks Coordination of Superior Level Staff Improvement(CAPES) – Finance Code 001. O.A. is partially supported by the Centre for Mathematics of the University of Coimbra (funded by the Portuguese Government through FCT/MCTES, DOI 10.54499/UIDB/00324/2020, DOI 10.54499/UIDP/00324/2020).
	P.G. thanks  Funda\c c\~ao para a Ci\^encia e Tecnologia FCT/Portugal for financial support through the projects UIDB/04459/2020 and UIDP/04459/2020. The authors would like to thank the referee for helping us to increase the reading of our article. This project has received funding from the European Research Council (ERC) under  the European Union's Horizon 2020 research and innovative programme (grant agreement   n. 715734). 
	
	\bibliographystyle{alea3}
	\bibliography{ALEA-AGS}
\end{document}